\documentclass[a4paper,11pt,makeidx]{amsart}

\usepackage{amssymb,eucal,graphicx}

\usepackage{ifthen}
\usepackage{epsfig}

\usepackage{enumerate}
\usepackage{psfrag}

\usepackage{amsmath}
\usepackage{amsfonts}
\usepackage{fancyhdr}
\usepackage[T1]{fontenc}

\RequirePackage{amsthm}
\RequirePackage{amssymb}
\usepackage[all,ps]{xy}
\usepackage{latexsym}

\numberwithin{equation}{section}

\newtheorem{theorem}[equation]{Theorem}
\newtheorem{thm}[equation]{Theorem}

\newtheorem{corollary}[equation]{Corollary}
\newtheorem{cor}[equation]{Corollary}

\newtheorem{lemma}[equation]{Lemma}

\newtheorem{proposition}[equation]{Proposition}
\newtheorem{prop}[equation]{Proposition}
\newtheorem{exa}[equation]{Example}
\newtheorem{remark}[equation]{Remark}

\def\Sing{\operatorname{Sing}}

\def\codim{\operatorname{codim}}

\def\im{\operatorname{im}}

\def\Exc{\operatorname{Exc}}

\def\Grass{\operatorname{Grass}}

\def\Sym{\operatorname{Sym}}
\def\slope{\operatorname{slope}}
\def\NE{\operatorname{NE}}
\def\x{\times}                   
\def\iso{\simeq}
\def\eqv{\equiv}
\def\sub{\subseteq}

\def\sup{\supseteq}

\def\+{\oplus}                   
\def\*{\otimes}                  
\def\hpil{\longrightarrow}       
\def\khpil{\rightarrow}

\def\mult{\operatorname{mult}}
\def\kod{\operatorname{kod}}

\def\Sym{\operatorname{Sym}}

\def\im{\operatorname{im}}

\def\Hom{\operatorname{Hom}}

\def\Pic{\operatorname{Pic}}

\def\Supp{\operatorname{Supp}}

\newcommand\Ii{{\mathcal I}}
\newcommand\Oc{{\mathcal O}}
\newcommand{\cO}{\mathcal{O}}

\newcommand\B{\mathcal B}
\newcommand\M{\mathcal M}
\newcommand\C{\mathcal C}

\newcommand\W{\mathcal W}
\newcommand\Ha{\mathcal H}

\newcommand\T{\mathcal T}
\newcommand\V{\mathcal V}
\newcommand\I{\mathcal I}

\newcommand\g{\mathfrak g}
\newcommand\ZZ{\mathbb{Z}}

\newcommand\NN{\mathbb{N}}
\newcommand\QQ{\mathbb{Q}}
\newcommand\RR{\mathbb{R}}

\newcommand\PP{\mathbb P}
\newcommand\Pp{\mathbb P}

\renewcommand{\P}{\mathbb{P}}
\renewcommand{\ss}{S^{[2]}}

\begin{document}

\title[On families of rational curves in the Hilbert square of a surface]{On families of rational curves in the Hilbert square of a surface
\\ \rm{(with an Appendix by Edoardo Sernesi)}}
\author[F.~Flamini, A.~L.~Knutsen, G.~Pacienza]{Flaminio Flamini*, Andreas Leopold Knutsen** and Gianluca Pacienza***}



\thanks{(*) and (***) Member of MIUR-GNSAGA at INdAM "F. Severi"}

\thanks{(**) Research supported by a Marie Curie Intra-European Fellowship within the 6th European
Community Framework Programme}

\thanks{(***) During the last part of the work the author benefitted from an "accueil en d\'el\'egation au CNRS"}

%

\vskip -30 pt

\begin{abstract} Under natural hypotheses we give an upper bound on the 
dimension of families of singular curves with hyperelliptic normalizations on a surface 
$S$ with $p_g>0$ via the study of the associated families of rational curves in $\ss$. 
We use this result to prove the existence 
of nodal curves of geometric genus $3$ with hyperelliptic normalizations, on a general $K3$ surface, 
thus obtaining specific $2$-dimensional families of rational curves in $\ss$.
We give two infinite series of examples 
of general, primitively polarized $K3$s such that their Hilbert squares contain a $\PP^2$ or a 
threefold birational to a $\PP^1$-bundle over a $K3$. 
We discuss the consequences on the Mori cone of the Hilbert square.
\end{abstract}

\maketitle

%
\section {Introduction} \label{S:intro}
%

For any smooth surface $S$, the Hilbert scheme $S^{[n]}$
parametrizing $0$-dimensional length $n$ subschemes of $S$ is a
smooth $2n$-dimensional variety whose inner geometry is naturally
related to that of $S$. For instance, if $\Delta \subset S^{[n]}$ is the
exceptional divisor, that is, the exceptional locus of the
Hilbert-Chow morphism $\mu: S^{[n]}\to \Sym^n(S)$, then 
irreducible (possibly singular) rational curves not contained in $\Delta$
roughly correspond to 
irreducible (possibly singular) curves on $S$ with a $\mathfrak{g}^1_{n'}$ on their normalizations, 
for some $n' \leq n$ (see \S\;\ref{S:corr} for the precise correspondence when $n=2$). 
One of the features of this paper is to show how ideas and techniques from one of the two sides of the correspondence
makes it possible 
to shed light on problems 
naturally arising on the other side.

If $S$ is a $K3$ surface, $S^{[n]}$ is a {\it hyperk\"ahler
manifold} (cf. \cite[2.2]{H1}) and rational curves play a fundamental r\^ole in the study of the (birational) 
geometry of $S^{[n]}$.  Indeed a result due to Huybrechts and Boucksom 
\cite{H2,Bou} implies in particular that these curves govern the ample cone 
of $S^{[n]}$ (we will recall the precise statement below and in \S\, \ref{S:prel}). 
The presence of 
a $\PP^n\subset S^{[n]}$ gives rise to a birational map 
(the so-called {\em Mukai flop} \cite{Muk}) to another hyperk\"ahler manifold and, for $n=2$, 
all birational maps between hyperk\"ahler fourfolds factor through a sequence of Mukai flops 
\cite{BHL,HY,W2,WW}. Moreover, as shown by 
Huybrechts \cite{H2}, uniruled divisors allow to describe the birational
K\"ahler cone of $S^{[n]}$ (see \S\, \ref{S:esempi} for the precise statement). For 
hyperk\"ahler fourfolds 
precise numerical and geometric properties of the rational 
curves that are extremal in the Mori cone have been conjectured by Hassett and Tschinkel \cite{HT}. 

The scope of this paper, and the structure of it as well, is twofold: we first devise general methods and tools to study families of curves with hyperelliptic normalizations on a surface $S$, mostly under the additional hypothesis that $p_g(S) >0$, in \S\;\ref{S:BB}-\S\;\ref{S:bound}. Then we apply these to obtain concrete results in the case of $K3$ surfaces, in \S\;\ref{S:esistenza3}-\S\;\ref{S:esempi}.
In particular, we have tried to develop a systematic
way to produce rational curves on $S^{[2]}$ by showing the existence
of nodal curves on $S$ with hyperelliptic normalizations.
 
To give an overview of the paper, we choose to start with the second part.

Let $(S,H)$ be a general, smooth, primitively polarized $K3$
surface of genus $p=p_a(H) \geq 2$. We have 
$N_1(\ss)_{\RR} \iso \RR[Y] \+ \RR[\PP^1_{\Delta}]$, where $\PP^1_{\Delta}$ is 
the class  of a rational curve in the ruling of
the exceptional
divisor $\Delta \subset \ss$, and $Y:= \{\xi\in \ss | \Supp(\xi)=\{p_0, y\}, \; \mbox{with}
\; p_0 \in S \; \mbox{and} \; y \in C \in |H|\}$, where $p_0$ and $C$ are chosen. One has that
$ \PP^1_{\Delta}$ lies on the boundary of the Mori cone and by the result of Huybrechts and Boucksom 
\cite{H2,Bou} mentioned above, if the Mori cone is closed, then also the other boundary is generated by the class of a rational curve. If $X \sim_{alg} aY -b \PP^1_{\Delta}$ is an irreducible curve in $\ss$, 
different from a  fiber of $\Delta$, then we define $a/b$ to be the {\it slope} of the curve. 
Thus, the lower the slope is, the closer is $X$ to the boundary of the Mori cone.
Describing
the Mori cone $\NE (\ss)$ amounts to computing
\[
 \slope(\NE (\ss)):= \inf \Big \{ \slope(X) \; | \; X \; \mbox{is an irreducible curve in} \; \ss \Big \},\]
and, if the Mori cone is closed, then $\slope(\NE (\ss))=\slope_{rat}(\NE (\ss))$, where 
\[
 \slope_{rat}(\NE (\ss)):= \inf \Big \{ \slope(X) \; | \; X \; \mbox{is an irreducible} \; rational \; \mbox{curve in} \; \ss \Big \}.
\]
(See \S\;\ref{S:prel}, \ref{S:class} and \ref{S:HEP} for further details.)

If now $C \in |mH|$ is an irreducible curve of geometric genus $p_g(C) \geq 2$ 
and with hyperelliptic normalization, let $g_0(C) \geq p_g(C)$ be the arithmetic genus of the {\it minimal} partial desingularization of $C$ that carries the $\g^1_2$
(see \S\;\ref{S:corr} and \S\;\ref{S:class}). By the unicity of the $\g^1_2$, $C$ defines a unique irreducible rational curve $R_C \subset \ss$ with class
 $R_C \sim_{alg} mY - ( \frac{g_0(C)+1}{2} ) \PP^1_{\Delta}$, cf. \eqref{eq:classR}.
(This formula is also valid if $R_C$ is associated to a given $\g^1_2$ on the normalization of an irreducible rational or elliptic curve $C$.) Thus, the higher $g_0(C)$ (or $p_g(C)$) is, and the lower $m$ is, the lower is the slope of $R_C$. This motivates the search for curves on $S$ with hyperelliptic normalizations of high geometric genus, thus ``unexpected'' from Brill-Noether theory.

It is well-known that there exist finitely many (nodal) rational curves, a one-parameter family 
of (nodal) elliptic curves, and a two-dimensional family of (nodal) curves of geometric genus $2$ in $|H|$ 
(see \S\;\ref{S:esistenza3}). Every such family yields in a natural way a two-dimensional family of irreducible rational curves in $\ss$, cf. \S\;\ref{S:BB}. Also note that, by a result of Ran \cite{ran},
the {\it expected dimension} of a family of rational curves in a symplectic fourfold, whence a posteriori also of a family of curves  with hyperelliptic normalizations lying on a $K3$, equals two 
(cf. Lemma \ref{cor:expdim}). In \cite[Examples 2.8 and 2.10]{fkp} we found 
two-dimensional families of nodal curves of geometric genus $3$ in $|H|$ having hyperelliptic 
normalizations when $p_a(H)=4$ or $5$. In this paper we generalize this:

\vspace{0,5cm}

\noindent {\bf Theorem \ref{thm:esistenza}.}
{\it Let $(S,H)$ be a general, smooth, primitively polarized $K3$
surface of genus $p=p_a(H) \geq 4$. Then the family of nodal curves in $|H|$
of geometric genus $3$ with hyperelliptic normalizations is nonempty, and each of its irreducible components is two-dimensional.}

\vspace{0,5cm}

The proof takes the whole \S\;\ref{S:esistenza3} and relies on a general principle
of constructing curves with hyperelliptic normalizations on general $K3$s outlined in
Proposition \ref{prop:degarg}:
first construct a marked $K3$ surface $(S_0,H_0)$ of genus $p$ such that $|H_0|$ contains a family
of dimension $\leq 2$ of nodal (possibly reducible) curves with the property that a desingularization
of some $\delta >0$ of the nodes is a limit of a hyperelliptic curve in 
the moduli space $\overline{\M}_{p-\delta}$ of stable curves of genus $p - \delta$ 
and such that this family is not contained in a higher-dimensional such family.
Then consider the parameter space $\W_{p,\delta}$ of pairs $((S,H),C)$, where $(S,H)$ is a smooth, primitively marked $K3$
surface of genus $p$ and $C \in |H|$ is a nodal curve with at least $\delta$ nodes. Now map (the 
local branches of) 
$\W_{p,\delta}$ into $\overline{\M}_{p-\delta}$ by partially normalizing the curves at $\delta$ of the nodes and 
mapping them to their respective classes. The existence of the particular family
in $|H_0|$ ensures that the image of this map intersects the hyperelliptic locus
$\overline{\mathcal{H}}_{p-\delta} \subset \overline{\M}_{p-\delta}$. A dimension count then shows that the dimension of the parameter space $\I \subset \W_{p,\delta}$ consisting of $((S,H),C)$ such that
a desingularization
of some $\delta >0$ of the nodes of $C$ is a limit of a hyperelliptic curve is at least $21$. Now the dominance on the $19$-dimensional moduli space of primitively marked $K3$ surfaces of genus $p$ follows
as the dimension of the special family on $S_0$ was $\leq 2$. 

The technical difficulties in the proof of Proposition \ref{prop:degarg} mostly arise because the 
curves in the special family on $S_0$ may be reducible (in fact, as in all arguments by degeneration, in practical applications they will very often be). Therefore we need to partially desingularize families of nodal curves, and this 
tool is provided in Appendix A by E.~Sernesi. Moreover, we need a careful study of the Severi varieties of {\it reducible} nodal curves on $K3$s, and here we use results of Tannenbaum \cite{Tan}. 

Given Proposition \ref{prop:degarg}, the proof of Theorem \ref{thm:esistenza} is then accomplished by constructing a suitable $(S_0,H_0)$ in Proposition \ref{prop:k3speciale} with $|H_0|$ containing a desired two-dimensional family of special curves, with $\delta=p-3$, and then showing that
the curves in the special family on $S_0$ in fact deform to curves with {\it precisely} $\delta$ nodes
on the general $S$ in Lemma \ref{lemma:k3speciale}. As will be discussed below, showing that the special family on $S_0$ is not contained in a family of higher dimension of curves with the same property, is quite delicate.

We also show that the associated rational curves in $\ss$ cover a threefold, cf. Corollary 
\ref{cor2:g=3}, and that $g_0=p_g=3$, cf. Remark \ref{rem:g=g0}.
Turning back to the description of $\NE (\ss)$, this shows that the class
of the associated rational curves in $\ss$ is 
$Y -  \frac{3}{2}  \PP^1_{\Delta}$, so that we obtain (cf. Corollary \ref{cor:g=3}): 
\vspace{0,5cm}
\[
\mbox{\eqref{eq:g=3boundinf}
\; \; \; \;  \; \;   \; \;   \; \;   \; \;  \; \;   \; \;   \; \;  \; \;  \; \;   \; \;   \;      
$\slope_{rat}(\NE (\ss)) \leq \frac {1}{2}$. \; \;  \; \;  \; \;  \; \;  \; \;   \; \;  \; \;  \; \;  \; \;  \; \; \; \;  \; \; }    
\]

\vspace{0,5cm}

In Propositions \ref{prop:IP2} and \ref{prop:IP1bundle} we present two infinite series of examples of {\it general}
primitively polarized $K3$ surfaces $(S,H)$ of infinitely many degrees
such that $\ss$ contains either a $\PP^2$ (these examples were shown to us by B.~Hassett) 
or a threefold birational to a $\PP^1$-bundle over a $K3$ and find the two-dimensional 
families of curves with hyperelliptic normalizations in $|H|$ corresponding to the lines and the fibres respectively. 
In particular, these examples show that the bound \eqref{eq:g=3boundinf}
can be improved for infinitely many degrees of the polarization. Namely, for any $n \geq 6$ and
$d \geq 2$, we get:
\vspace{0,5cm}
\[
\mbox{\eqref{eq:slope-n}
\; \; \; \;  \; \;   \; \;   \; \;   \; \;  \; \;     
$\slope_{rat}(\NE (\ss)) \leq \frac {2}{2n-9}$ \; if \;$p = p_a(H) =
n^2-9n+20$; \; \;  \; \;  \; \;  \; \;  \; \;   \; \;  \; \;  \; \;  \; \;  \; \; \; \;  \; \; }    
\]
\[
\mbox{\eqref{eq:slope-d}
\; \; \; \;  \; \;   \; \;   \; \;   \; \;  \; \;     
$\slope_{rat}(\NE (\ss)) \leq \frac {1}{d}$ \; if \;$p = p_a(H) =d^2$. \; \;  \; \;  \; \;  \; \;  \; \;   \; \;  \; \;  \; \;  \; \;  \; \; \; \;  \; \; }    
\]

\vspace{0,5cm}
\noindent Nevertheless, to our knowledge, \eqref{eq:g=3boundinf} is  the first non-trivial bound valid for any genus 
$p$ of the polarization.  

The proofs of Propositions \ref{prop:IP2} and 
\ref{prop:IP1bundle} are again by deformation, but unlike the proof of Proposition \ref{prop:degarg}, we now deform $\ss_0$ of a special $K3$ surface $S_0$. The idea is to start with a special quartic surface $S_0 \subset \PP^3$ such that
$\ss_0$ contains a $\PP^2$ or a threefold birational to a $\PP^1$-bundle over itself, perform the standard involution on $\ss_0$ to produce a new such and then deform $\ss_0$ keeping the new one by keeping a suitable polarization on the surface that is different from $\cO_{S_0}(1)$. Here we use 
results on deformations of symplectic fourfolds by Hassett and Tschinkel \cite{HT} and Voisin \cite{V2}.

By a result proved in \cite{fkp}, any irreducible curve $C \in |H|$ with hyperelliptic normalization must satisfy 
$g_0(C) \leq \frac{p+2}{2}$, where $p = p_a(H) $ 
(cf. Theorem \ref{thm:ijm} and \eqref{rho<0}). It is then natural to ask whether this inequality actually ensures the existence of 
such curves. We call this ``The hyperelliptic existence problem'' and we see that a positive solution to this 
problem would yield a bound on the slope of rational curves that is much stronger than the ones obtained above, 
cf. \eqref{eq:boundinf}. In this sense, Theorem \ref{thm:esistenza} is hopefully only the first step towards 
stronger existence results.

The study of curves on $S$ with hyperelliptic normalizations is not the only way to obtain bounds on the slope of the Mori cone of $\ss$. In fact, an irreducible curve $C \in |mH|$ with a singular point $x$ of multiplicity
$\mult_x(C)$ yields an irreducible curve in $\ss$ with class $mY - (1/2)\mult_x(C)$ (see the proof
of Theorem \ref{thm:eff}). In particular, if $p = p_a(H)$, one has the bound (cf. Theorem \ref{thm:seshadribound})
\vspace{0,5cm}
\[
\mbox{ \eqref{eq:seshadribound}
\; \; \; \;  \; \;   \; \;   \; \;   \; \;  \; \;   \; \;   \; \;  \; \;  \; \;   \; \;   \;      
$\slope(\NE (\ss)) \leq \sqrt{\frac{2}{p-1}}$, \; \;  \; \;  \; \;  \; \;  \; \;   \; \;  \; \;  \;
 \;  \; \;  \; \; \; \;  \; \; }    
\]

\vspace{0,5cm}
\noindent obtained by using well-known results on Seshadri constants on $S$. This bound  is stronger than
\eqref{eq:g=3boundinf} but weaker than the bounds on the slope of the Mori cone obtained from
\eqref{eq:slope-n} and 
\eqref{eq:slope-d}. Moreover, one relatively easily sees that the best bound one can obtain by 
Seshadri constants is in any case weaker than \eqref{eq:slope-n} and 
\eqref{eq:slope-d} and also weaker than the ones one could obtain by solving ``The hyperelliptic existence problem'', cf. \eqref{eq:boundinf}. In any case, note that \eqref{eq:seshadribound}, \eqref{eq:slope-n} and 
\eqref{eq:slope-d} show that the bounds tend to zero as the degree of the polarization tends to infinity, that is,  
\[
\mbox{ \eqref{eq:inf0}
\; \; \; \;  \; \;   \; \;   \; \;   \; \;  \;       
$\inf \Big\{ \slope(\NE (\ss)) \; | \; S$ \; is a projective $K3$ surface $\Big\}=0$, 
 \; \;  \; \;  \; \;  \; \;  \; \;   \; \;  \; \;  \; \;  \; \;  \; \; \; \;  \; \; }    
\]
and likewise for $\slope_{rat}(\NE (\ss))$.

All the families of curves in $|H|$ with hyperelliptic normalizations we have seen above have in fact dimension equal to two, the expected one. Moreover, a crucial point in the proof of Theorem \ref{thm:esistenza}
is to bound the dimensions of families of irreducible curves with hyperelliptic normalizations on the special $K3$ surface $S_0$. This brings us over to the description of the first part of this paper.

The problem of bounding the dimension of special families of curves on surfaces, like in our case of curves with hyperelliptic normalizations, is interesting in its own, may be 
studied for larger classes of surfaces, and may lead to further applications in other contexts. 
Whereas  methods from adjunction theory have proved very useful for the study of {\it smooth}
hyperelliptic curves on surfaces \cite{se,sv,BFL}, these methods do not extend to the case of {\it singular} curves, where in fact very little seems to be known. Even in the relevant case of nodal curves on smooth surfaces, whose parameter spaces (the so-called {\it Severi varieties}) have received much attention over the years and have been studied also in relation with moduli problems (see e.g. \cite{Ser1} for 
$\Pp^2$ and \cite{F2} for surfaces of general type), 
the dimension of their subloci consisting of curves with hyperelliptic normalizations
is not determined.  

The precise question we address is whether there exists an upper bound on the dimension of families of irreducible curves on a projective surface with hyperelliptic normalizations. One easily sees that, if the canonical system of the surface is birational,
 then no curve with hyperelliptic normalization can move, cf. e.g. \cite{kn}. On the other 
hand, taking 
any surface $S$ admitting a (generically) $2:1$ map onto  a rational surface $R$ and pulling back the families of rational curves on $R$, we obtain families of arbitrarily high dimensions of curves on $S$ having hyperelliptic 
normalizations. Moreover, the infinite series of examples in Proposition \ref{prop:IP2} of general,
primitively polarized $K3$ surfaces $(S,H)$ such that $\ss$ contains a $\PP^2$ shows that one cannot even hope, in general, 
to find a bound in the simplest case of Picard number one: in fact, the $(3m-1)$-dimensional family of rational curves
in $|\cO_{\PP^2}(m)|$ yields a $(3m-1)$-dimensional family of irreducible curves in $|mH|$ having hyperelliptic normalizations,
cf. \S\;\ref{S:hassett}. Nevertheless, for a large class of surfaces, it is possible
to derive a geometric consequence on the family $V$, when its dimension
is greater than two:

\vspace{0,5cm}

\noindent {\bf Theorem \ref{thm:apppg}'.}
{\it Let $S$ be a smooth, projective surface with $p_g(S) >0$. Let $V$ be a reduced and irreducible scheme parametrizing a flat family of irreducible curves on $S$ with hyperelliptic normalizations (of genus $\geq 2$) such that $\dim (V) \geq 3$.
Then the algebraic equivalence class $[C]$ of the curves parametrized by $V$ has a decomposition $[C] = [D_1] + [D_2]$ into algebraically moving classes
such that $[D_1+D_2] \in \overline{V}$. Moreover 
the rational curves in $S^{[2]}$ corresponding to the irreducible curves parametrized by $V$ cover only a (rational) surface $R \subset S^{[2]}$.}

\vspace{0,5cm}

In fact we prove a stronger result, cf. 
Theorem \ref{thm:apppg}, that in particular relates the decomposition $[C] = [D_1] + [D_2]$ to the $\g^1_2$s on 
the normalizations of the curves parametrized by $V$. This additional point will in fact be the crucial one in 
our application in the proof of Theorem \ref{thm:esistenza}. An immediate corollary is that 
the ``na{\"i}ve'' dimension bound one may hope for, thinking about the fact that rational curves in $\ss$ arising from curves on $S$ of geometric genera $\leq 2$ move in dimension at most two, is in fact
true under additional hypotheses on $V$, cf. Corollary \ref{cor:spezza}. These are satisfied if e.g. the 
N\'eron-Severi group of $S$ is of rank $1$ and generated by the class of a curve in $V$, and seem quite natural, taking into account the examples of large families mentioned above. 

 The idea of the proof of Theorem \ref{thm:apppg} is rather simple and geometric and  illustrates
well the rich interplay between the properties of curves on $S$
and those of subvarieties of $\ss$. The proof relies on the following two fundamental results:

The first is Mori's bend-and-break technique (see Lemma \ref{lemma:BBiper} for the precise version we need), which gives a breaking into reducible members of a family of rational curves of dimension $\geq 3$ covering a surface.

The second is a suitable version of Mumford's well-known theorem on
$0$-cycles on surfaces with $p_g>0$ (cf. Corollaries \ref{cor:mumford} and \ref{cor:mumford2}). The consequence of particular interest to us is that any threefold in $\ss$
can only carry a two-dimensional covering family of rational curves when $p_g(S) >0$, cf. 
Proposition \ref{cl:dimRV}.

Combining those two ingredients, we see that any family satisfying the hypotheses of Theorem \ref{thm:apppg} yields a family
of rational curves in $\ss$ of the same dimension $\geq 3$, that can therefore only cover a surface in $\ss$, on which we can apply bend-and-break to produce a reducible member. Then we have to show that we can also produce a decomposition
of the curves on $S$ into algebraically moving classes, and this is carried out in Proposition \ref{prop:algdec}.

Beside the application in the proof of Theorem \ref{thm:esistenza}, we hope that Theorem \ref{thm:apppg} and the ideas behind its
 proof will find more applications. One is a Reider-like
result for families of singular curves with hyperelliptic normalizations obtained in \cite{kn}, where also more examples are given.

The paper is organized as follows. We go from the more general results to those peculiar to the case of $K3$ surfaces. We start in \S\,\ref{S:BB} with the correspondence between curves with hyperelliptic
normalizations on any smooth surface $S$ and  rational curves on $\ss$ and prove other preliminary results, before turning to the bend-and-break lemma for families of rational curves covering a surface in $\ss$.
The version of Mumford's theorem we need for our purposes is proved in \S\,\ref{S:mumford}, and then rephrased in terms of rational quotients. Then we prove (a stronger version of) Theorem \ref{thm:apppg}' in \S\,\ref{S:bound}. We then turn to 
$K3$ surfaces and prove Theorem \ref{thm:esistenza} along the lines of the degeneration argument sketched above. 
Section \ref {S:K3}, apart from some 
known facts on the Hilbert scheme of points on a $K3$ surface, contains the computation of the classes of rational curves 
in $\ss$
associated to curves in $S$ with rational, elliptic or hyperelliptic normalizations, as explained in \S\;\ref{S:corr}. 
The relation 
between the existence of such a curve, its singular Brill-Noether number (an invariant introduced in \cite{fkp}) and 
the slope of the Mori cone of 
$\ss$ is also discussed, as well as the relation between the slope of the Mori cone and Seshadri constants. 
We end the paper presenting the two series of examples of general $K3$ surfaces whose Hilbert square contains a $\P^2$
(respectively a threefold birational to a $\P^1$-bundle over a $K3$) and discussing the 
numerical properties of  a line (respectively a fibre) in it, as well as those of the associated singular curves in $S$ with 
hyperelliptic normalizations.
In Appendix \ref{A:edoardo} by E.~Sernesi the reader will find a general result 
about partial desingularizations of families of nodal curves. 

\vspace{0,5cm}

{\it Acknowledgements.} The authors thank L.~Caporaso, O.~Debarre, A.~Iliev and A. Verra for  
useful discussions related to these problems. We are extremely grateful to: C.~Ciliberto, for many valuable
conversations and helpful comments on the subject and for having pointed out some 
mistakes in a preliminary version of this paper; B.~Hassett, for having 
pointed out the examples behind Proposition \ref{prop:IP2};
E.~Sernesi, for many helpful conversations and for 
his Appendix \ref{A:edoardo}. 
We finally express our gratitude to the Department of Mathematics, 
Universit\`a "Roma Tre" and to the Institut de Recherche Math\'ematique Avanc\'ee, Universit\'e L. Pasteur et CNRS,
where parts of this work have been done,
for the nice and warm atmosphere as well as for the
kind hospitality.

%
\section {Rational curves in $\ss$}\label{S:BB}
%

Let $S$ be a smooth, projective surface.
In this section we gather some basic results that will be needed in the rest of the paper. We first 
describe the natural correspondence between rational curves in $\ss$ and curves on $S$ with 
rational, elliptic or hyperelliptic normalizations. Then, 
in \S\;\ref{SS:BB}, we apply Mori's bend-and-break technique to rational curves in $\Sym^2(S)$ covering a surface. 

Recall that we have the natural Hilbert-Chow morphism 
$\mu: \ss \to \Sym^2(S)$ that resolves $\Sing (\Sym^2(S)) \iso S$. The $\mu$-exceptional divisor $\Delta \subset \ss$ is a $\PP^1$-bundle over $S$. The Hilbert-Chow morphism gives an obvious one-to-one correspondence between irreducible curves in $\ss$ not contained in $\Delta$ and irreducible curves in $\Sym^2(S)$ not contained in $\Sing (\Sym^2(S))$. We will therefore often switch back and forth between working on $\ss$ and $\Sym^2(S)$.

%
\subsection{Irreducible rational curves in $\ss$ and curves on $S$} \label{S:corr}
%

Let $T \subset S \x \ss$ be the incidence variety, with projections $p_2: T \to \ss$ and $p_S: T \to S$.
Then $p_2$ is finite of degree two, branched along $\Delta \subset \ss$. (In particular, $T$ is smooth as
$\Delta$ is.)

Let  $X \subset \ss$ be an irreducible rational curve not contained in $\Delta$.
We will now see how $X$ is equivalent to one of three sets of data on $S$.

Let 
$\nu_X: \widetilde{X} \iso \PP^1 \to X$ be the normalization and set $X':=p_2^{-1}(X) \subset T$. By the universal property of blowing up, we obtain a commutative square:
\begin{equation} \label{eq:diarami}
\xymatrix{
\widetilde{C}_X \ar[r]^f \ar[d]_{\tilde{\nu}_X} & \widetilde{X} \ar[d]^{\nu_X} & \hspace{-1cm} \iso \PP^1 \\
X' \ar[r]^{{p_2}_{|X'}} & X, 
}
\end{equation}
defining the curve $\widetilde{C}_X$, $\tilde{\nu}_X$ and $f$. In particular,
$\tilde{\nu}_X$ is birational and $\widetilde{C}_X$ admits a $\g^1_2$ (i.e., a $2:1$ morphism onto $\PP^1$, given by $f$), but may be singular, or even reducible.
Set $\tilde{\nu}:={p_S}_{|X'} \circ \tilde{\nu}_X: \widetilde{C}_X \to S$.

Assume first that $\widetilde{C}_X$ is irreducible.

We set $C_X:= \tilde{\nu}(\widetilde{C}_X) \subset S$. Since
$X \not \subset \Delta$, $C_X$ is a curve. As $\widetilde{C}_X$ carries  a $\g^1_2$, it is easily seen that also the normalization of $C_X$ does, that is, $C_X$ 
has rational, elliptic or hyperelliptic normalization. Moreover, it is easily seen that $\tilde{\nu}: \widetilde{C}_X \to C_X$ is generically of degree one. Indeed, for general $x \in C_X$, as $x \not \in p_S(p_2^{-1}(\Delta))$, we can write
$({p_S}_{|X'}) ^{-1}(x)= \{ (x, x+y_1), \ldots, (x, x+y_n)\}$, where $n:=\deg \tilde{\nu}$.
By definition of $p_2$, and since  $X'=p_2^{-1}(X)$, we must have
that each $(y_i,x+y_i) \in X'$, for $i=1, \ldots, n$, and each couple 
$((x,x+y_i), (y_i,x+y_i))$ is the pushdown by $\tilde{\nu}_X$ of an element of the $\g^1_2$ on 
$\widetilde{C}_X$. Hence, each couple $(x,y_i)$ is the pushdown by the normalization morphism 
of an element of the induced $\g^1_2$ on the normalization of 
$C_X$. Since $x$ has been chosen general, $x \not \in \Sing (C_X)$, so that we must have $n=1$, as claimed.

In particular, by construction, $\tilde{\nu}: \widetilde{C}_X \to C_X$ is a partial desingularization
of $C_X$, in fact, it is the {\it minimal} partial desingularization of $C_X$ carrying the $\g^1_2$
in question (which is unique, if $p_g(C_X) \geq 2$). 

We have therefore obtained:

\begin{itemize}
\item[(I)]  the data of  an irreducible curve $C_X \subset  S$ together with 
a partial normalization $\tilde{\nu}: \widetilde{C}_X \to C_X$ with a $\g^1_2$ on
$\widetilde{C}_X$ (unique, if  $p_g(C_X) \geq 2$), such that $\tilde{\nu}$ is minimal with respect to
the existence of the $\g^1_2$.
\end{itemize}

Next we treat the case where $\widetilde{C}_X$ is reducible. 
In this case, it must consist of two irreducible smooth rational components, 
$\widetilde{C}_X= \widetilde{C}_{X,1} \cup \widetilde{C}_{X,2}$, that are identified by $f$.

If $\tilde{\nu}$ does not contract any of the components, set 
$C_{X,i}:= \tilde{\nu}(\widetilde{C}_{X,i}) \subset S$ and
$n_{X,i}:= \deg \tilde{\nu}_{|\widetilde{C}_{X,i}}$, for $i=1,2$. We therefore obtain:

\begin{itemize}
\item[(II)] the data of a curve $C_X= n_{X,1}C_{X,1} + n_{X,2}C_{X,2} \subset  S$, with
 $n_{X,i} \in \NN$, $C_{X,i}$ an irreducible, rational curve, a morphism
$\tilde{\nu}: \widetilde{C}_X=  \widetilde{C}_{X,1} \cup \widetilde{C}_{X,2}
\to C_{X,1} \cup C_{X,2}$ (resp. $\tilde{\nu}: \widetilde{C}_X \to C_{X,1}$ if $C_{X,1}=C_{X,2}$)
that is $n_{X,i}:1$ on each component and where $\widetilde{C}_{X,i}$ is the normalization of $C_{X,i}$,
and an identification morphism
$f: \widetilde{C}_{X,1} \cup \widetilde{C}_{X,2} \iso \PP^1 \cup \PP^1 \to \PP^1$.
\end{itemize}

If $\tilde{\nu}$ contracts one of the two components of $\widetilde{C}_X$, say
$\widetilde{C}_{X,2}$, to a point $x_X \in S$ (it is easily seen that it cannot contract both), then
$\mu(X) \subset \Sym^2(S)$ is of the form $\{x_X+C_X \}$, for an irreducible curve 
$C_X \subset  S$, which is necessarily rational. It is easily seen that 
$C_X= \tilde{\nu}(\widetilde{C}_{X,1})$ and 
$\deg \tilde{\nu}_{|\widetilde{C}_{X,1}}=1$, so that we obtain:

\begin{itemize}
\item[(III)] the data of an irreducible rational curve $C_X \subset S$ together with a point $x_X \in S$.
\end{itemize}

Note that in all cases (I)-(III), the support of the curve $C_X$ on $S$ is simply
\begin{equation}\label{eq:aiuto0}
 \Supp (C_X)= \mbox{one-dimensional part of} \; \{x \in S \; | \; x \in \Supp (\xi) \; \mbox{for some} \;
\xi \in X \}
\end{equation}
and the set is already purely one-dimensional precisely unless we are in case (III) with
$x_X \not \in  C$. 

Conversely from the data (I), (II) or (III) one recovers an irreducible rational curve
in $\ss$ not contained in $\Delta$.
Indeed, in case (I) (resp. (II)), the $\g^1_2$ on $\widetilde{C}_X$ (respectively, the identification
$f$) induces a $\PP^1 \subset \Sym^2(\widetilde{C}_X)$ and this is mapped to an irreducible rational curve in $\Sym^2(S)$ by the natural composed morphism 
\[
\xymatrix{
\Sym^2(\widetilde{C}_X) \ar[r]^{\tilde{\nu}^{(2)}} & \Sym^2(C_X)  \ar@{^{(}->}[r] & \Sym^2(S).
}
\]
The irreducible rational curve $X \subset \ss$ is the strict transform by $\mu$ of this curve.
In case (III), $X \subset \ss$ is the strict transform by $\mu$ of $\{x_X+C_X\} \subset \Sym^2(S)$.

We see that the data (III) correspond precisely to rational curves of type 
$\{x_0+C\} \subset \Sym^2(S)$, where $x_0 \in S$ is a point and $C \subset S$ is an irreducible rational curve.
Moreover, it is easily seen that the data (II) correspond precisely to the images by
\[ \alpha: \widetilde{C}_1 \x \widetilde{C}_2 \iso \PP^1 \x \PP^1 \hpil C_1+C_2 \subset \Sym^2(S), \]
resp.
\[ \alpha: \Sym^2(\widetilde{C}) \iso \PP^2 \hpil \Sym^2(C) \subset \Sym^2(S), \]
of irreducible rational curves in $|n_1F_1+n_2F_2|$ for $n_1, n_2 \in \NN$, resp.
$|nF|$ for an integer $n \geq 2$, where $\Pic (\widetilde{C}_1 \x \widetilde{C}_2) \iso \ZZ[F_1] \+ \ZZ[F_2]$,
resp. $\Pic (\Sym^2(\widetilde{C})) \iso \ZZ[F]$, and $C_1, C_2$, resp. $C$, are irreducible rational curves on $S$ and ``\;$\widetilde{}$\;'' denotes normalizations. The data of type (II) will however not 
be studied more in this paper, where we will focus on the other two, mostly on (I).

Note that an irreducible rational curve $X \subset \Sym^2(S)$ arising from rational (resp. elliptic) curves $C$ as in case (I) moves in
$\Sym^2(C)$, which is a surface birational to $\PP^2$ (resp. an elliptic ruled
surface), and a curve $X \subset \Sym^2(S)$ of the form $\{x_X+C\}$ moves in the threefold 
$\{S + C\}$, which is birational to a $\PP^1$-bundle over $S$, and contains $\Sym^2(C)$.

At the same time, it is well-known that if $\kod(S) \geq 0$, then rational 
curves on $S$ do not move and elliptic curves move in at most one-dimensional families. 
This follows for instance from the following general result (that we will later need in the case $p_g=2$). 

\begin{lemma} \label{lemma:dimre}
 Let $S$ be a smooth, projective surface with $\kod (S) \geq 0$ containing an $n$-dimensional irreducible family of irreducible curves of geometric genus $p_g$. Then $n \leq p_g$ and
if equality occurs, then either 
the family consists of a single smooth rational 
curve; or $\kod (S) \leq 1$ and $n \leq 1$; or 
$\kod (S) =0$.
\end{lemma}

\begin{proof}
This is ``folklore''. For a proof see \cite{kn}.
\end{proof}

As a consequence, if $\kod(S) \geq 0$, then rational curves in $\Sym^2(S)$ arising from rational or elliptic curves on $S$ move in families of dimension at most two in $\Sym^2(S)$.

On the
other hand, irreducible rational curves $X \subset \Sym^2(S)$ arising from curves on $S$ with hyperelliptic normalizations of geometric genus $p_g \geq 2$ (necessarily of type (I)),  
move in a family
whose dimension equals that of the family of curves with hyperelliptic
normalizations in which $C \subset S$ moves 
(by unicity of the $\mathfrak{g}^1_2$). Apart from some special cases, it is easy to see that the converse is also true:

\begin{lemma}\label{lem:famcost}
Let $\{X_b\}_{b\in B}$ be a one-dimensional irreducible family of irreducible rational curves in $\Sym^2(S)$
covering a (dense subset of a) proper, reduced and irreducible  surface $Y \subset \Sym^2(S)$
that does not coincide with $\Sing(\Sym^2(S)) \cong S$.

Then  $C=C_{X_b}$ in $S$ for every  $b \in B$ (notation as above) if and only if either
$Y= \Sym^2(C_0)$, with either $C_0 \subset S$ an irreducible rational curve and $C \eqv nC_0$
for $n \geq 1$, or $C_0=C \subset S$ an irreducible elliptic curve; or 
$Y= C + C':= \{ p+p' \; |\;  p\in C,\; p'\in C'\}$, 
with $C$ an irreducible rational curve and $C' \subset S$ any irreducible curve; or
$Y=C_1+C_2$, with $C_1, C_2 \subset S$ irreducible rational curves and
$C=n_1C_1+n_2C_2$ for $n_1, n_2 \in \NN$.
\end{lemma}

\begin{proof} The "if" part is immediate. For the converse, we treat the three cases (I)-(III) separately. 

If $C$ is as in (I), then clearly $Y \subset \Sym^2(C)$, so that $Y = \Sym^2(C)$ and $C$ must be either 
rational or elliptic, as $Y$ is uniruled. 

If $C = n_1C_1 + n_2 C_2$ as in (II), then either $C_1 = C_2 = : C_0$ and again 
$Y = \Sym^2(C_0)$, or $C_1 \neq C_2$ and $Y = C_1 + C_2$. 

Finally, if $C$ is as in (III), then, for every $b \in B$, we have $\{X_b\}_{b \in B} = \{x_b + C\}_{b\in B}$ for 
some $x_b \in S$, and the $\{x_b\}_{b \in B}$ define the desired curve $C'$.  

\end{proof}

We note that by Lemma \ref{lemma:dimre} also the rational curves in $\Sym^2(S)$ arising from
singular curves of geometric genus $2$ on $S$ move in at most two-dimensional families. We will see below that this is a general phenomenon, under some additional hypotheses.
We will focus our attention on curves with hyperelliptic normalizations (of genus $p_g \geq 2$) in 
Sections \ref{S:bound}-\ref{S:esempi}.

%
\subsection{Bend-and-break in $\Sym^2(S)$}\label{SS:BB} 
%

Let $V \sub \Hom (\PP^1,\Sym^2(S))$ be a reduced and irreducible subscheme 
(not necessarily complete). We consider the universal map
\begin{equation} \label{eq:univ2}
\xymatrix{ \mathcal{P}_V := \PP^1 \x V  \ar[r]^{ \hspace{0,5cm}\Phi_V} & \Sym^2(S) }
\end{equation} 
and assume that the following two conditions hold:
\begin{equation} \label{eq:juve}
\mbox{For any} \; v \in V, \; \Phi_V(\PP^1 \x v ) \not \sub \Sing(\Sym^2(S)) \iso S; \; \mbox{and}
\end{equation} 
\begin{equation} \label{eq:roma}
\Phi_V \; \mbox{is generically finite}
\end{equation} 
(the latter just means that $V$ induces a flat family of rational curves in $\Sym^2(S)$ of dimension 
$\dim(V)$).
Set
\begin{equation} \label{eq:defrv}
R_V:= \overline{\im (\Phi_V)},
\end{equation}
the
Zariski closure of $\im (\Phi_V)$ in $\Sym^2(S)$. It is the (irreducible) uniruled subvariety of $\Sym^2(S)$ covered by the curves parametrized by $V$. 
In the language of \cite[Def. 2.3]{Kol}, $R_V$ is the closure of the {\it locus} of the family 
$\Phi_V$.
Note that $\dim (R_V) \geq 2$ if $\dim (V) \geq 1$ by \eqref{eq:roma}. Moreover (cf. e.g. \cite[Prop. 2.1]{HT1}),
\begin{equation} \label{eq:dimkod0}
\dim (R_V) \leq 3 \; \mbox{if} \; \kod(S) \geq 0.
\end{equation}

When $R_V$ is a surface, using Mori's  bend-and-break technique we obtain the following result. 
In the statement we underline the fact that the breaking can be made in such a way that,
 for general $\xi, \eta \in R_V$, two components of the reducible (or non-reduced) member at the border of the family pass through
$\xi$ and $\eta$, respectively. This will be central in our applications 
(Proposition \ref{prop:algdec} 
and \S\;\ref{S:esistenza3}, where we prove Theorem \ref {thm:esistenza}). We give the proof because we could not find in literature precisely the 
statement we will need.

\begin{lemma} \label{lemma:BBiper}
 Assume that $\dim (V) \geq 3$ and $\dim(R_V)=2$.

Let $\xi$ and $\eta$ be any two distinct general points of $R_V$.
 Then there is a curve $Y_{\xi,\eta}$ in $R_V$ such that
 $Y_{\xi,\eta}$ is algebraically equivalent to $(\Phi_V)_*(\PP^1_v)$ and
either
\begin{itemize}
\item[(a)] there is an irreducible nonreduced component of $Y_{\xi,\eta}$
containg $\xi$ and $\eta$; or
\item[(b)] there are two distinct, irreducible components of
$Y_{\xi,\eta}$
           containg $\xi$ and $\eta$, respectively.
\end{itemize}
\end{lemma}

\begin{proof}
Since $\dim (V) \geq 3$ by assumption, by \eqref{eq:roma} we can  pick a
one-dimensional smooth subscheme $B=B_{\xi,\eta} \subset V$ parametrizing
curves in $V$ such that $(\Phi_V)_* (\PP^1 \x v)$ contains both $\xi$
and $\eta$, for every $v \in B$.  We therefore
have a family of rational
curves:
\begin{equation} \label{eq:famrazpassanti}
\Phi_B :=(\Phi_V){|_B}: \PP^1 \x B \hpil R_V.
\end{equation}
and two marked (distinct) points $x,y \in \PP^1$ such that $\Phi_B(x \x
B)=\xi$ and $\Phi_B(y \x B) = \eta$, such that each $\Phi_B(\PP^1
\x v)$ is nonconstant, for any $v \in B$; in particular $\Phi_B(
\PP^1 \x B)$ is a surface.

As in the proofs of \cite[Lemma 1.9]{km} and \cite[Cor. II.5.5]{Kol}, let 
$\overline{B}$ be any smooth compactification of $B$. 
Consider the surface $\PP^1 \x \overline{B}$. Let $0 \in \overline
B$ denote a point at the boundary, $\P^1_0$ the fibre over $0$ of
the projection onto the second factor and $x_0, y_0 \in \P^1_0
\subset \PP^1 \x \overline{B}$ the corresponding marked points. 
By the {\it Rigidity Lemma} \cite[Lemma 1.6]{km},
$\Phi_B$ cannot be defined at the point $x_0$, as in the proof of \cite[Cor. 1.7]{km},
and the same argument
works for $y_0$. 

Therefore, to resolve the indeterminacies of the rational map
$\Phi_B: \PP^1 \x \overline{B} --\khpil R_V$, we must
at least blow up  $\PP^1 \x \overline{B}$ at the points $x_0$ and
$y_0$. Now let $W$ be the blow-up of $\PP^1 \x \overline B$ such
that $\overline{\Phi}_B: W \hpil R_V$ is an extension
of $\Phi_B$, that is, we have a commutative diagram
\[
\xymatrix{
    W   \ar@{->}[d]_{\pi}      \ar@{->}[dr]^{\overline{\Phi}_B} &   \\
                        \PP^1 \x \overline{B} \ar@{-->}[r]^{\Phi_B} &  R_V.
}
\]
Let $E_{x_0}:= \pi^{-1}(x_0)$ and $E_{y_0}:= \pi^{-1}(y_0)$. Note that neither of these can be
contracted by $\overline{\Phi}_B$, for otherwise $\Phi_B$ itself would be defined at
$x_0$ or $y_0$.

Therefore the curve $\overline{\Phi}_B(E_{x_0})$ has an irreducible component $\Gamma_{\xi}$ containing $\xi$
and the curve $\overline{\Phi}_B(E_{y_0})$ has an irreducible component $\Gamma_{\eta}$ containing $\eta$
and by construction, $\Gamma_{\xi}+\Gamma_{\eta} \sub {\overline{\Phi}_B}_*(\pi^{-1}(\PP^1 \x 0))$ and the latter
is the desired curve $Y_{\xi,\eta}$. The two cases (a) and (b) occur as
$\Gamma_{\xi}=\Gamma_{\eta}$ or $\Gamma_{\xi} \neq \Gamma_{\eta}$, respectively.
\end{proof}

%
\section {Rationally equivalent zero-cycles on surfaces with $p_g>0$} \label{S:mumford}
%

In this section we extend to the singular case a consequence of Mumford's
result on zero-cycles on surfaces with $p_g>0$ (cf.
\cite[Corollary p.~203]{Mum}) and reformulate the results in terms of rational quotients.

%
\subsection{Mumford's Theorem}\label{SS:mum}
%

The main result of this subsection, which we prove in detail for the reader's convenience, 
relies on the following generalization of Mumford's result (cf. \cite[Chapitre 22]{V}  and references therein,
for a detailed account).

\begin{thm}\label{mumford} (see \cite[Prop. 22.24]{V})
Let $T$ and $Y$ be smooth projective varieties. Let
$Z\subset Y\times T$ be a cycle of codimension equal to
$\dim (T)$. Suppose there exists a subvariety $T'\subset T$ of dimension $k_0$
such that, for all $y\in Y$, the
zero-cycle $Z_y$ is rationally equivalent in $T$ to a cycle supported on $T'$.

Then, for all $k>k_0$ and for all $\eta\in H^0 (T, \Omega^k_T)$, we have
$$[Z]^*\eta =0 \; {\rm in} \; H^0 (Y, \Omega^k_Y)$$where, as costumary,
$[Z]^* \eta$ denotes the differential form
induced on $Y$ by the correspondence $Z$.
\end{thm}

Mumford's original ``symplectic''
argument and the theorem above yield the following result
(see \cite[Corollary p.~203]{Mum}).

\begin{corollary}\label{cor:mumford}
Let $S$ be a smooth, irreducible projective surface with $p_g(S)> 0$ and 
$\Sigma \subset S^{[n]}$ a reduced, irreducible (possibly
singular) complete subscheme such that $\mu(\Sigma) \not \subset \Sing(\Sym^n(S))$, where
$\mu: S^{[n]} \to \Sym^n(S)$ is
the {\em Hilbert-Chow morphism}.

If there exists a subvariety $\Gamma \subset \Sym^n(S)$ such
that $\dim(\Gamma) \leq 1$, $\Gamma \not \subset
\Sing(\Sym^n(S))$ and all the zero-cycles parametrized
by $\mu(\Sigma)$ are rationally equivalent to zero-cycles
supported on $\Gamma$, then $\dim(\Sigma) \leq n$.
\end{corollary}

\begin{proof}
Let $\pi : \widetilde{\Sigma} \to \Sigma \subset S^{[n]}$ be the
desingularization morphism of $\Sigma$. Let $Z = \Lambda_{\pi}
\subset \widetilde{\Sigma} \times S^{[n]}$ be the graph of $\pi$.
Then $Z \cong \widetilde{\Sigma}$, so that $\codim(Z) = \dim(S^{[n]})$,
as in Theorem \ref{mumford}.
By assumption, $\mu(\Sigma)$ parametrizes zero-cycles of length
$n$  on $S$ that are all rationally equivalent to zero-cycles
supported on $\Gamma$, with $\dim(\Gamma) \leq 1$.
Since $\mu(\Sigma)$ is not contained in $\Sing(\Sym^n(S))$
 by assumption,
$\mu|_{\Sigma}: \Sigma \to \mu(\Sigma)$ is birational. If
$\Gamma'$ denotes the strict transform of $\Gamma$ under $\mu$, we
get that $\dim(\Gamma') \leq 1$.

We can apply Theorem \ref{mumford} with $Z = Y = \widetilde{\Sigma}$,
$T = S^{[n]}$ and $T' = \Gamma'$. Thus, for each $k >1$ and for
each $\eta \in H^0(\Omega^k_{S^{[n]}})$, $[Z]^*\eta = 0$ in
$H^0( \widetilde{\Sigma}, \Omega^k_{\widetilde{\Sigma}})$.

Let $\omega \in H^0(S, K_S)$ be a non-zero $2$-form on $S$. As in
\cite[Corollary]{Mum}, we define:
\[
\omega^{(n)} := \sum_{i=1}^n
p_i^*(\omega) \in H^0 (S^n , \Omega^2_{S^n})
\]
where $S^n$ is the
$n^{th}$-cartesian product and $p_i$ is the natural projection
onto the $i^{th}$ factor, $1 \leq i \leq n$. The form
$\omega^{(n)}$ is $\Sym(n)$-invariant and, since we have that
$\mu$ is surjective, this induces a canonical $2$-form $\omega^{[n]}_{\mu}
\in H^0(S^{[n]}, \Omega^2_{S^{[n]}} )$ (see \cite[\S 1]{Mum}, where 
$\omega^{[n]}_{\mu} = \eta_{\mu}$ in the notation therein). From what we observed above, $[Z]^*(\omega^{[n]}_{\mu}) = 0$ as a
form in $H^0(\widetilde{\Sigma}, \Omega^2_{\widetilde{\Sigma}})$.
Consider
$$
 (\Sym^n (S))_0 := \Big \{ \xi = \sum_{i=1}^n x_i \; | \; x_i \neq x_j,
 \; 1 \leq i \neq j \leq n\ {\rm and\ such\ that}\
 \omega(x_i)\in \Omega^2_{S,x_i}\ {\rm is\ not}\ 0 \Big\}.
$$
Then $(\Sym^n (S))_0 \subset \Sym^n (S)$ is an open dense
subscheme that is isomorphic to its preimage via $\mu $ in
$S^{[n]}$. For each $\xi \in (\Sym^n(S))_0$, $\xi$ is a smooth
point and$$\pi_n: S^n \to \Sym^n (S)$$is {\'e}tale over $\xi$.
Thus, the $2$-form $\omega^{(n)} \in H^0(S^n, \Omega^2_{S^n})$ is
non-degenerate on the open subset $(S^n)_0$ of points in the
preimage of $(\Sym^n(S))_0$, i.e. it defines a non-degenerate
skew-symmetric form on the tangent space of $(S^n)_0$.

Let $\pi^0_n := \pi_n|_{(S^n)_0}$; since $\pi^0_n: (S^n)_0 \to
(\Sym^n (S))_0 $ is  {\'e}tale,  there exists a
$2$-form$$\omega_0^{(n)} \in H^0((\Sym^n (S))_0,
\Omega^2_{(\Sym^n(S))_0})$$such that $\omega^{(n)} =
\pi_n^*(\omega_0^{(n)}) $ and $\omega_0^{(n)}$ is also
non-degenerate. Therefore, the maximal isotropic subspaces of
$\omega_0^{(n)}(\xi)$ are $n$-dimensional.

Now $\Sigma \subset S^{[n]}$ and $\Sigma \cap \mu^{-1}((\Sym ^n
(S))_0) \neq \emptyset$, since $\mu(\Sigma) \not \subset 
\Sing(\Sym^n (S))$ by assumption.
Since $\Sigma$ is reduced, let $\xi \in \Sigma \cap \mu^{-1}
((\Sym^n (S))_0) $ be a smooth point. Then, since $\Sigma_{smooth} =
\pi^{-1}(\Sigma_{smooth})$, by abuse of notation we still denote 
by $\xi  \in \widetilde{\Sigma}$ the corresponding point.
We know that $[Z]^* \omega^{[n]}_{\mu}(\xi) = 0$ in the tangent space
$T_{\xi}(\widetilde{\Sigma})$. Since
$$
\xi \in \Sigma_{smooth} \cap
\mu^{-1}((\Sym^n (S))_0) \subset (\Sym^n (S))_0,
$$
then
$[Z]^*(\omega^{[n]}_{\mu})  = \omega_0^{(n)}|_{\Sigma_{smooth} \cap
\mu^{-1}((\Sym^n (S))_0)}$. This implies $\dim(\Sigma) \leq n$.
\end{proof}

%
\subsection{The property $RCC$ and rational quotients}\label{SS:RCC}
%

Recall that a variety $T$ (not necessarily proper or smooth) is said to be
{\it rationally chain connected}  ($RCC$, for brevity), if
for each pair of very general points $t_1, t_2 \in T$ there exists a connected
curve $\Lambda \subset T$ such that $t_1, t_2 \in \Lambda$ and
each irreducible component of $\Lambda$ is rational (see \cite{Kol}). Furthermore, by
\cite[Remark 4.21(2)]{Deb}, if $T$ is proper and $RCC$,
then {\it each pair} of points can be joined by a
connected chain of rational curves.

Also recall that, for any smooth variety $T$, there exists a variety $Q$, called the {\it rational quotient of} $T$, together with a
rational map 
\begin{equation} \label{eq:f}
f:T - - \to Q,
\end{equation}
whose very general fibres are
equivalence classes under the $RCC$-equivalence relation
(see, for instance, \cite[Theorem 5.13]{Deb} or \cite[IV, Thm. 5.4]{Kol}).

In this language, an equivalent statement 
of Corollary \ref{cor:mumford} is:

\begin{corollary} \label{cor:mumford2} 
Let $S$ be a smooth, projective surface with $p_g(S) >0$. If $Y \subset S^{[n]}$ is a 
complete subvariety of dimension $>n$ not contained in $\Exc (\mu)$,
then any 
desingularization of $Y$ has a rational quotient of dimension at least two.
\end{corollary}

\begin{proof} 
Let $\widetilde{Y}$ be any desingularization of $Y$ and $Q$ its rational quotient. Up to
resolving the indeterminacies of $f:\widetilde{Y}- - \to Q$,
we may assume that $f$ is a
proper morphism whose very general fibre is a $RCC$-equivalence
class, so that in particular {\it each} fibre is $RCC$
(see \cite[Thm. 3.5.3]{Kol}).

If $\dim(Q) = 0$, it follows that $\widetilde{Y}$ (so also
$Y$) is $RCC$, contradicting Corollary \ref{cor:mumford}.

If $\dim(Q)=1$, then by cutting $\widetilde{Y}$ with
$\dim (Y)-1$ general very ample divisors,
we get a curve  $\Gamma'$  that intersects every fibre of
$f$. Every point of $\widetilde{Y}$ is
connected by a chain of rational curves to some point on
$\Gamma'$. We thus obtain a contradiction by Corollary \ref{cor:mumford}
(with $\Gamma$ the image of $\Gamma'$ in $\Sym^2(S)$).
\end{proof}

Let now $R_V$ be the variety covered by a family  of rational curves in $\Sym^2(S)$ parametrized by $V$, as 
defined in \eqref{eq:defrv},  
$\widetilde{R}_V$ be any desingularization of $R_V$ and 
$Q_V$ be the rational quotient of $\widetilde{R}_V$. 
Of course $\dim (Q_V) \leq \dim (R_V)-1$, as $R_V$ is uniruled by construction. 

\begin{lemma} \label{lem:QV}
If $\dim(V) \geq \dim (R_V)$, then $\dim (Q_V) \leq \dim (R_V) -2$ (for any desingularization $\widetilde{R}_V$ of $R_V$).
In particular, if $\dim(V) \geq 2$ and $\dim(R_V) = 2$, then any
desingularization of $R_V$ is a rational surface.
\end{lemma}

\begin{proof}
With notation as in \S\;\ref{SS:BB}, we have $\dim({\mathcal P}_V) \geq \dim (R_V)+1$, so that
the general fibre of $\Phi_V$ is at least one-dimensional, cf. \eqref{eq:univ2}. This means
that, if $\xi$ is a general point of $R_V$, there exists a family
of rational
curves in $R_V$ passing through $\xi$, of dimension $\geq 1$.
Of course the same is true for a
general point of $\widetilde{R}_V$. Thus,
the very general fibre of $f$ in \eqref{eq:f}
has dimension at least two, whence $\dim(Q_V) \leq \dim (R_V)-2$.
The last statement  follows from the fact that any smooth surface that is
$RCC$ is rational (cf. \cite[IV.3.3.5]{Kol}).
\end{proof}

Combining Corollary \ref{cor:mumford2} and Lemma \ref{lem:QV}, we then get:

\begin{prop} \label{cl:dimRV}
If $p_g(S) >0$ and $\dim (V) \geq 2$, then either
\begin{itemize}
\item[(i)] $R_V$ is a surface with rational desingularization; or
\item[(ii)] $\dim (V)=2$, $R_V$ is a threefold and any desingularization of $R_V$ has a two-dimensional rational quotient.
\end{itemize}
\end{prop}

\begin{proof}
By \eqref{eq:dimkod0}, $\dim (R_V)=2$ or $3$.
If $\dim (R_V)=2$, then (i) holds by Lemma \ref{lem:QV}.
If $\dim (R_V)=3$, then $\dim (Q_V) = 2$ by Corollary \ref{cor:mumford2}. Hence $\dim (V)=2$ by 
Lemma \ref{lem:QV} and (ii) holds. 
\end{proof}

\begin{remark} \label{rem:apppg2} {\rm 
Let $S$ be a smooth, projective surface with $p_g(S) >0$ and let $Y \subset S^{[2]}$ be a uniruled threefold different from $\Exc(\mu)$, where $\mu: S^{[2]} \to \Sym^2(S)$ is the Hilbert-Chow morphism. 

Take a covering family $\{C_v\}_{v \in V}$ of rational curves on $Y$. By Corollary \ref{cor:mumford2} the family must be two-dimensional (see Lemma \ref{lem:QV}).
Then the curves in the covering family yield, via the correspondence described in 
\S\;\ref{S:corr}, curves on $S$ with rational, elliptic or hyperelliptic normalizations, and the correspondence is one-to-one in the hyperelliptic case. We therefore see 
that we must be in one of the following cases:
\begin{itemize}
\item[(a)] $S$ contains an irreducible rational curve $\Gamma$ and
\[ Y=\{ \xi \in \ss \; | \; \mbox{$\Supp(\xi) \cap \Gamma \neq \emptyset$} \};
\]
\item[(b)] $S$ contains a one-dimensional irreducible family $\{ E \}_{v \in V}$ of irreducible
elliptic curves and 
\[
Y = \overline{\{ \xi \in E_v^{[2]} \}}_{v \in V }; \] or
\item[(c)] $S$ contains a two-dimensional, irreducible family of irreducible
curves with hyperelliptic normalizations, not contained in a higher dimensional irreducible family,
and $Y$ is the locus covered by the corresponding rational curves in $\ss$.  
\end{itemize}(Note that in fact case (b) can only occur for $\kod(S) \leq 1$ by Lemma \ref{lemma:dimre} and case (c) only when $|K_S|$ is not birational. The latter fact is easy to see, cf. e.g. \cite{kn}.) 

In
the case of $K3$ surfaces, uniruled divisors play a particularly
important r\^ole \cite[\S 5]{H2}, cf. \S\;\ref{S:esempi}. Now all cases (a)-(c) above
occur on a general, projective $K3$ surface with a polarization of genus $\geq 6$. 
In fact, cases (a) and (b) occur on any projective $K3$
surface since it necessarily contains a one-dimensional family of irreducible, elliptic curves and a zero-dimensional family of rational curves, by a well-known theorem of Mumford
(see the proof in \cite[pp.~351-352]{MM} or \cite[pp.~365-367]{BPV}). Case (c) occurs on a general primitively polarized $K3$ surface of genus $p \geq 6$ by Corollary \ref{cor2:g=3} below
with a family of curves of geometric genus $3$. In addition to this, in Proposition \ref{prop:IP1bundle} we will see that there is another threefold as in (c) arising from curves of geometric genus $>3$ in the hyperplane linear system on general projective $K3$
surfaces of infinitely many degrees. 

Moreover, there is not a one-to-one correspondence between families as in (a), (b) or (c) above and uniruled threefolds in
$S^{[2]}$. In fact, in Proposition \ref{prop:IP2} we will see that there is a two-dimensional family of 
curves with hyperelliptic normalizations, as in (c), in the hyperplane linear systems on general $K3$ surfaces of infinitely 
many degrees whose associated rational curves cover only a $\PP^2$ in $S^{[2]}$.}
\end{remark}

%
\section {Families of curves with hyperelliptic normalizations} \label{S:bound}
%

The purpose of this section is to study the dimension of families of curves on 
a smooth projective surface $S$  with hyperelliptic normalizations. 

We first remark that it is not difficult to see that if $|K_S|$ is birational, 
then the dimension of such a family is forced to be zero (see e.g. \cite{kn}). At the same time it is easy to find 
obvious examples of surfaces, even with $p_g(S) >0$, with large families of
curves with hyperelliptic normalizations, namely surfaces admitting a finite $2:1$ map
onto a rational surface. (For examples of such cases, see e.g. 
\cite{Hor1, Hor2, Hor3,Hor4,SD,se,sv,BFL} to mention a few.)
In these cases one can pull back the families of rational curves on the rational surface to obtain families of 
curves on $S$ with hyperelliptic normalizations of arbitrarily high dimensions. Moreover, in Proposition \ref{prop:IP2} 
below we will see that even a general, primitively polarized $K3$ surface $(S,H)$, for infinitely many degrees, 
contains a $\PP^2$ in its Hilbert square, which is  not contained in $\Delta$ (but the surface is not a double cover of a 
$\PP^2$, by generality).   Therefore, by the correspondence in \S\;\ref{S:corr},  
$S$ contains large families of curves with hyperelliptic normalizations. 
One can see that in all these examples of large families the algebraic equivalence class of the members 
breaks into nontrivial effective decompositions. For example, in the mentioned $K3$ case of Proposition \ref{prop:IP2}, 
we will see that the  curves in $|\mathcal{O}_{\PP^2}(n)|$ in $\PP^2 \subset \ss$ correspond to curves in $|nH|$. In this section we will see that this is a general phenomenon, with the help of Lemma \ref{lemma:BBiper}.

To this end, let  $V$ be a reduced and irreducible scheme parametrizing a flat family of curves on $S$ all having constant geometric genus
$p_g \geq 2$ and hyperelliptic normalizations. 
Let $\varphi:\C \khpil V$ be the universal family. 
Normalizing $\C$ we obtain, possibly restricting to an open dense subscheme
of $V$, a flat family $\tilde{\varphi} : \widetilde{\C}\to V$ of smooth hyperelliptic 
curves of genus $p_g \geq 2$ (cf. \cite[Thm. 1.3.2]{tes}). 
Let $\omega_{\widetilde{\C}/V}$ be the relative dualizing sheaf. 
As in \cite[Thm. 5.5 (iv)]{LK},  consider the morphism 
$\gamma: \widetilde{\C}\to \PP(\tilde{\varphi}_* (\omega_{\widetilde{\C}/V}))$ over $V$. This morphism 
is finite and of relative degree two onto its image, which we denote by $\mathcal{P}_V$.
We thus obtain a universal family $\psi: \mathcal{P}_V \khpil V$ 
of rational curves mapping  to  $\Sym^2(S)$, as in 
\eqref{eq:univ2}, satisfying \eqref{eq:juve} and \eqref{eq:roma}.  (Strictly speaking,
\eqref{eq:univ2} denoted a universal family of {\it maps}, whereas it now denotes a universal family of
{\it curves}.) 
To summarize, recalling \eqref{eq:defrv}, we have
\begin{equation} \label{eq:univ}
\xymatrix{ & \widetilde{\C} \ar[dl]_{\pi} \ar[dr]_{\tilde{\varphi}} \ar[r]^{\gamma}  & 
\mathcal{P}_V \ar[d]^{\psi} \ar[r]^{\Phi_V} & R_V \\
S  & & V. &
}
\end{equation}

Also note that \eqref{eq:univ} is compatible with the correspondence of case (I) in \S \;\ref{S:corr}, in the sense
that, for general $v \in V$, we have (using the same notation as in \S \;\ref{S:corr})
\begin{equation} \label{eq:univa}
 \pi( \tilde{\varphi}^{-1}(v))= p_S(p_2^{-1}X_v)=(p_S)_*(p_2^{-1}X_v) =C_{X_v}, \; \mbox{with} \; 
X_v= \mu^{-1}_* \Big(\Phi_V(\psi^{-1}(v))\Big) \subset \ss,
\end{equation}
where $\mu$ is the Hilbert-Chow morphism (in particular, $p_S$ and $p_2$ are the first and second projections, 
respectively, from the incidence variety $T \subset S \x \ss$). Note that the second equality in \eqref{eq:univa} 
follows as $p_S$ is generically one-to-one on the curves in question, as we saw in \S \;\ref{S:corr}. 
This will be central in the proof of the next result.

We now apply Lemma \ref{lemma:BBiper} to ``break'' the curves on $S$.

\begin{prop} \label{prop:algdec}
Let $S$ be a smooth, projective surface and $V$ 
and $R_V$ as above. Assume that $\dim (V) \geq 3$ and $\dim(R_V)=2$ and let $[C]$ be the algebraic equivalence class 
of the members 
parametrized by $V$.

Then there is a decomposition into two effective, algebraically moving classes
\[
[C] = [D_1] + [D_2]
\]
such that, for general $\xi, \eta \in R_V$,
there are
effective divisors $D'_1  \sim_{alg} D_1$ and $D'_2 \sim_{alg}
D_2$ such that $\xi \subset D'_1$ and $\eta \subset D'_2$ and $[D'_1 +
D'_2] \in \overline{V}$, 
where $\overline {V}$ is the closure of $V$ in the component of the Hilbert scheme
of $S$ containing $V$.
\end{prop}

\begin{proof} 
For general $\xi, \eta \in R_V$, both being supported at two distinct points on $S$, 
let $B=B_{\xi,\eta} \subset{V}$ be as in the proof of Lemma \ref{lemma:BBiper}
and  $\overline{B}$ be any smooth compactification of $B$. By abuse of notation, we will consider $\xi$ and $\eta$ as being points in $\ss$. By (the proof of) Lemma \ref{lemma:BBiper},
using the Hilbert-Chow morphism, there is a flat family $\{X_b \}_{b \in \overline{B}}$ of curves in 
the surface $\mu^{-1}_*(R_V) \subset \ss$ (where $\mu$ is the Hilbert-Chow morphism as usual)
parametrized by $\overline{B}$, such that, for general $b \in B$, $X_b$ is an irreducible rational curve and
\begin{equation} \label{eq:pd}
C_{X_b}= (p_S)_*(p_2^{-1}(X_b)) = \pi( \tilde{\varphi}^{-1}(b)),
\end{equation}
with notation as in \S\;\ref{S:corr} (cf. \eqref{eq:univa}). In particular, $\{C_{X_b}\}_{b \in B}$ is a one-dimensional nontrivial subfamily of 
the family $\{C_{X_v}\}_{v \in V}$ given by $V$. Moreover, for some $b_0 \in \overline{B} \setminus B$, we have
$X_{b_0} \sup Y_{\xi} + Y_{\eta}$,
where $Y_{\xi}$ and $Y_{\eta}$ are irreducible rational curves (possibly coinciding) 
such that $\xi \in Y_{\xi}$ and $\eta \in Y_{\eta}$. Also note that
$Y_{\xi}, Y_{\eta} \not \subset \Delta \subset \ss$. 

Pulling back to the incidence variety $T \subset S \x \ss$, we obtain a flat family 
$\{X'_b:=p_2^{-1}(X_b) \}_{b \in \overline{B}}$ of curves in $T$, such that
\begin{equation} \label{eq:pd2}
 X'_{b_0}:= p_2^{-1}(X_b) \sup p_2^{-1}(Y_{\xi}) + p_2^{-1}(Y_{\eta}) =: Y_{\xi}'+ 
Y_{\eta}'. 
\end{equation}

Note that the family $\{X'_b\}_{b \in \overline{B}}$ is in fact a family of curves in the
incidence variety
$T_0 \subset S \x \mu^{-1}_*(R_V)$, which is a surface contained in $T$. Since
$p_S$ maps this family to a family of curves covering (an open dense subset of) $S$, by \eqref{eq:pd}, we see that
$(p_S)_{|T_0}$ is surjective, in particular generically finite.  
Thus, choosing $\xi$ and $\eta$ general enough, we can make sure they lie outside of the images by $p_2$ of the finitely many curves contracted by $(p_S)_{|T_0}$. Hence $q^{-1}(Y_{\xi})$ and 
$q^{-1}(Y_{\eta})$ are not contracted by $p_S$. 

Therefore, recalling \eqref{eq:pd} and \eqref{eq:pd2} and letting $b' \in B$ be a general point, we get
\[ C \sim_{alg} (p_S)_* X'_{b'}  \sim_{alg} (p_S)_* X'_{b_0} \sup (p_S)_* Y_{\xi}'+ 
(p_S)_* Y_{\eta}' \sup D_{\xi} +D_{\eta},\]
where 
$D_{\xi}:= p(q^{-1}Y_{\xi})$ and $D_{\eta}:= p(q^{-1}Y_{\eta})$. 

By construction we have $D_{\xi} \supset \xi $ and $D_{\eta} \supset \eta$,  viewing
$\xi$ and $\eta$ as length-two subschemes of $S$. (Note that $D_{\xi}$ and
$D_{\eta}$ are not necessarily distinct.) Possibly after adding additional components to
$D_{\xi}$ and $D_{\eta}$, we can in fact assume that  
\[ C \sim_{alg} (p_S)_* X'_{b'} = D_{\xi} +D_{\eta},\]
 with $D_{\xi}$ and $D_{\eta}$ not necessarily reduced and irreducible.
Since this construction can be repeated for general
$\xi, \eta \in R_V$ and the set $\{x \in S \; |\; \; x \in \Supp (\xi) \;\; \mbox{for some} \; \xi \in R_V \}$
is dense in $S$, as the curves parametrized by $V$ cover the whole surface $S$, the obtained curves
$D_{\xi}$ and $D_{\eta}$ must move in an algebraic system of dimension at least one. 

By construction, $D_{\xi} + D_{\eta}$ lies in the border of the family 
$\varphi: \C \khpil V$ of curves on $S$, and as such, $[D_{\xi} + D_{\eta}]$ lies in the closure of $V$ in the component of the Hilbert scheme
of $S$ containing $V$. 
Moreover, as the number of such decompositions is finite (as $S$ is projective and the divisors are effective), we can find one decomposition
$ [C] = [D_1] + [D_2]$ holding for general $\xi, \eta \in R_V$.  
\end{proof}

The next two results are immediate consequences:

\begin{thm} \label{thm:apppg}
 Let $S$ be a smooth, projective surface with $p_g(S) >0$. Then the following conditions are equivalent:
\begin{itemize}
\item[(i)]  $S^{[2]}$ contains an  irreducible   surface $R$ with rational desingularization,
such that  $R \neq \mu_*^{-1}(C_1 + C_2), \; \mu_*^{-1}(\Sym^2(C))$ for rational curves 
$C, C_1, C_2 \subset S$  and $R \not \subset \Exc(\mu)$, where $\mu: S^{[2]} \to \Sym^2(S)$ is the Hilbert-Chow morphism;
\item[(ii)] $S$ contains a flat family of irreducible curves with hyperelliptic normalizations of geometric genus  $p_g \geq 3$, parametrized by a reduced and irreducible scheme $V$  
such that $\dim (V) \geq 3$.
\end{itemize}
Furthermore, if any of the above conditions holds, then
\begin{itemize}
\item[(a)] the rational curves in $\ss$ that correspond to the irreducible curves parametrized by $V$, cover only the surface  $R$  in $S^{[2]}$; and
\item[(b)] the algebraic equivalence class $[C]$ of the curves parametrized by $V$ has an effective decomposition $[C] = [D_1] + [D_2]$ into algebraically moving classes
such that, for general $\xi, \eta \in R$, there are effective divisors $D'_1  \sim_{alg} D_1$
and $D'_2 \sim_{alg} D_2$ such that $\xi \subset D'_1$, $\eta \subset D'_2$ and $[D'_1+D'_2] \in \overline{V}$, where $\overline {V}$ is the closure of $V$ in the component of the Hilbert scheme
of $S$ containing $V$.
\end{itemize}
\end{thm}

\begin{proof}
Assume (ii) holds. By Proposition \ref{cl:dimRV} we have that $R_V \subset \Sym^2(S)$ is a  surface with rational desingularization, so that (i) holds.

Assume now that (i) holds. Then $R$ carries a family of rational curves of dimension
$n \geq 3$. By Lemma \ref{lem:famcost} and the assumptions in (i), this yields an $n$-dimensional family of curves on $S$ that have rational, elliptic or hyperelliptic normalizations.
From Lemma \ref{lemma:dimre}, we get (ii).

Finally, assume that these conditions hold. Then (a) follows from Proposition \ref{cl:dimRV} again, where $R$ is the proper transform via $\mu$ of the surface $R_V$ therein; finally, (b) follows from Proposition \ref{prop:algdec}.
\end{proof}

\begin{cor} \label{cor:spezza}
Let $S$ be a smooth, projective surface with $p_g(S)>0$ and $V$ be a
reduced, irreducible scheme parametrizing a flat family of
irreducible curves with hyperelliptic normalizations (of geometric genus $\geq 2$).
Denote by $[C]$ the algebraic equivalence class of the members of $V$.

If $[C]$ has no decomposition into effective, algebraically moving classes,  then
$\dim (V) \leq 2$.
\end{cor}
In particular, Corollary \ref{cor:spezza} holds 
when e.g. $NS(S) = \ZZ[C]$.

The examples with the double covers of smooth rational surfaces and the result in 
Proposition \ref{prop:IP2} mentioned above, show that the results above are natural.

The statement in Theorem \ref{thm:apppg}(b) shows that in fact the length-two zero-dimensional schemes on the curves 
in the family corresponding to the elements of the $\mathfrak{g}^1_2$s on their normalization, are in fact 
``generically cut out'' by moving divisors in a fixed algebraic decomposition of the class of the members in the family. This reminds of the nowadays well-known results of Reider and their generalizations
\cite{re,BFS,BS}. In fact, Theorem \ref{thm:apppg}(b) can be used to prove 
a Reider-like result involving the arithmetic and geometric genera 
of the curves in the family,  cf. \cite{kn}. 
Moreover, the precise statement in Theorem 
\ref{thm:apppg}(b) will be crucial in the next section, where we will prove existence of curves with hyperelliptic 
normalizations by degeneration methods.

%
\section{Nodal curves of geometric genus $3$ with hyperelliptic normalizations on $K3$ surfaces}\label{S:esistenza3}
%

In the rest of the paper we will focus on the existence of curves with ``Brill-Noether special'' hyperelliptic
normalizations (i.e. of geometric genera $>2$) and in this section we will  see that Theorem \ref{thm:apppg}(b)
is particularly suitable to prove existence results by degeneration arguments.

To do this and to discuss some consequences on $\ss$, we will in the rest of the paper focus on $K3$ surfaces,
which
in fact were one of our original motivations for this work.

We start with the following observation combining a result of Ran,
already mentioned in the Introduction, with the results from the previous section.

\begin{lemma} \label{cor:expdim}
Let $S$ be a smooth, projective
$K3$ surface and $L$ be a globally generated line bundle of sectional genus $p \geq 2$ on $S$.
Let $|L|^{hyper} \sub |L|$ be the subscheme parametrizing irreducible curves in $|L|$ with hyperelliptic normalizations.

Then, any irreducible component of $|L|^{hyper}$ has dimension
$\geq 2$, with equality holding if $L$ has no decomposition into moving classes.
\end{lemma}

\begin{proof}
Any $n$-dimensional component of $|L|^{hyper}$ yields an $n$-dimensional family of
irreducible rational curves in $\ss$. By \cite[Cor. 5.1]{ran}, we have $n \geq 2$.
The last statement follows from  Corollary \ref{cor:spezza}.
\end{proof}

\noindent

The main aim of this section is to apply Theorem \ref{thm:apppg}(b) to prove:

\begin{theorem}\label{thm:esistenza}
Let $(S,H)$ be a general, smooth, primitively polarized $K3$
surface of genus $p=p_a(H) \geq 4$. Then the family of nodal curves in $|H|$
of geometric genus $3$ with hyperelliptic normalizations is nonempty, and each of its irreducible components is two-dimensional.
\end{theorem}

In \cite{fkp} we studied which linear series may appear on normalizations of
irreducible curves on $K3$ surfaces. To do so, we introduced a singular Brill-Noether number 
$\rho_{sing} (p_a, r, d, p_g)$ whose negativity, 
when $\Pic(S) \iso \ZZ[H]$, ensures non-existence of curves in $|H|$, with $p_a = p_a(H)$ and of geometric genus $p_g$, 
having normalizations admitting a $\mathfrak{g}^r_d$ (we will return to this in \S\;\ref{S:HEP} below). 
Moreover, in \cite[Examples 2.8 and 2.10]{fkp}, we already gave examples of nodal curves with
hyperelliptic normalizations with geometric genus $3$ and arithmetic genus $4$ or $5$. 
Theorem \ref{thm:esistenza} shows  that this is a general phenomenon. 
The proof will  be given in the remainders of this section. Moreover,  
we will also determine the dimension of the
locus covered in $\ss$  by the rational curves associated to curves in a component of the family:

\begin{cor}  \label{cor2:g=3}
Let $(S,H)$ be a general, smooth, primitively polarized $K3$
surface of genus $p=p_a(H) \geq 6$. Then the subscheme of $|H|$ parametrizing
nodal curves of geometric genus $3$ with hyperelliptic normalizations contains a two-dimensional component 
$V$ such that $\dim(R_V)=3$.
\end{cor}

This corollary in particular shows that all three cases in Remark \ref{rem:apppg2} occur on a general $K3$
surface.
In \S\;\ref{S:class}-\ref{S:HEP} we will both compute the classes of the corresponding rational curves in $\ss$ 
(see \eqref{eq:genere3})
and discuss some of the consequences of Theorem \ref{thm:esistenza} on the Mori cone of $\ss$.

Before starting on the proof of Theorem \ref{thm:esistenza}, we
recall that, for any smooth surface $S$ and any line
bundle $L$ on $S$, such that $|L|$ contains smooth, irreducible
curves of genus $p:= p_a(L)$, and any positive integer $\delta \leq p$, one
denotes by $V_{|L|, \delta}$ the locally closed and
functorially defined subscheme of $|L|$ parametrizing the
universal family of irreducible curves in $|L|$ having
$\delta$ nodes as the only singularities and, consequently, geometric genus $p_g:= p- \delta$.
These are classically called {\em Severi varieties} of irreducible, $\delta$-nodal
curves on $S$ in $|L|$.

It is nowadays well-known, as a direct consequence of  Mumford's theorem on the existence of nodal
rational curves on $K3$ surfaces (see the proof in
\cite[pp.~351-352]{MM} or \cite[pp.~365-367]{BPV}) and standard
results on Severi varieties,  that if
$(S,H)$ is a general, primitively polarized $K3$ surface of genus
$p \geq 3$, then the Severi variety $V_{|H|, \delta}$ is nonempty
and {\it regular}, i.e. it is smooth and of the expected dimension
$p-\delta$, for each $\delta \leq p$ (cf. \cite[Lemma 2.4 and Theorem 2.6]{Tan}; see also
e.g. \cite{CS,F}).

The regularity property follows from the fact that, since by definition $V_{|L|, \delta}$ parametrizes irreducible curves,
the nodes of these curves impose independent conditions on $|L|$ (cf. \cite{CS,F} and 
\cite[Remark 2.7]{Tan}). From equisingular deformation theory, this implies that suitable obstructions to some
locally trivial deformations are zero. In other words, it implies first that, for any $\delta' > \delta$,
$V_{|L|, \delta'} \subset \overline{V}_{|L|, \delta}$ (see \cite[Anhang F]{Sev}, \cite{Wa} and \cite[Thm. 4.7.18]{Ser2} for $\PP^2$ and \cite[\S\,3]{Tan} for $K3$s). Furthermore,
if $[C] \in V_{|L|, \delta + k}$, $k >0$, is a general point of an irreducible component, the fact that the nodes
impose independent conditions allows to clearly describe what $\overline{V}_{|L|, \delta}$ looks like locally around
the point $[C]$: it is the union of $\delta+k \choose \delta$ smooth branches through $[C]$, each branch
corresponding to a choice of $\delta$ "marked" (or "assigned") nodes among the $\delta+k$ nodes of $C$, and
these branches intersect transversally at $[C]$; moreover, the other $k$ "unassigned" nodes of $C$ disappear
when one deforms $[C]$ in the corresponding branch of $\overline{V}_{|L|, \delta}$ (see \cite[Anhang F]{Sev}, 
\cite{Wa} and \cite[\S\,1]{Ser1} for $\Pp^2$ and \cite[\S\,3]{Tan} for $K3$s).

The situation is slightly different for {\em reducible}, nodal curves in $|L|$. Since they appear in the proof of Theorem
\ref{thm:esistenza}, we also have to take care of this case. To this end, we define the ``degenerated'' version of $V_{|L|, \delta}$ by
\begin{eqnarray}\label{eq:w}
  W_{|L|, \delta} & := & \Big \{ C \in |L| \; |  \; C, \mbox{ not necessarily irreducible,
  has only nodes } \\
\nonumber  &   &  \mbox{as singularities and {\em at least} $\delta$ nodes} \Big  \}.
\end{eqnarray}For the same reasons as above, $W_{|L|, \delta}$ is a
locally closed subscheme of $|L|$. Note that
\begin{equation} \label{eq:irred}
  W_{|L|, \delta} = \cup_{\delta' \geq \delta} V_{|L|, \delta'} \; \mbox{if all the curves in $|L|$ are irreducible},
\end{equation}which is a partial compactification of $V_{|L|, \delta}$.

Let $[C] \in W_{|L|, \delta}$. Choosing any subset 
$\{p_1, \ldots, p_{\delta}\}$ of $\delta$ of its nodes, one obtains a {\em pointed curve} 
$(C; p_1 , \ldots, p_{\delta})$, where $p_1, \ldots, p_{\delta}$ are also called the 
{\em marked} (or {\em assigned}) nodes of $C$ (cf. \cite[Definitions 3.1-(ii) and 3.6-(i)]{Tan}).

Recall that there exists an algebraic scheme, which we denote by
\begin{equation}\label{eq:branch}
\B(C; p_1, p_2 , \ldots, p_{\delta}),
\end{equation}locally closed in $|L|$, representing the functor of infinitesimal deformations of $C$ in $|L|$ that preserve the marked nodes, i.e. the functor of locally trivial infinitesimal deformations of the pointed
curve $(C; p_1 , \ldots, p_{\delta})$ (cf. \cite[Proposition 3.3]{Tan}, where we have identified 
the schemes therein with their projections into the linear system $|L|$). 
In other words, $\B(C; p_1, p_2 , \ldots, p_{\delta})$ is the {\em local branch} of
$W_{|L|, \delta}$ around $[C] \in W_{|L|, \delta}$, corresponding to the choice of the $\delta$ marked nodes. 
We have:

\begin{theorem}\label{thm:Tan3} (cf. \cite[Theorem 3.8]{Tan}) Let 
$(C; p_1 , \ldots, p_{\delta})$ be as above. Assume that the general element of $|L|$ is a smooth, irreducible 
curve and that the partial normalization of $C$ at the $\delta$ marked nodes $p_1 , \ldots, p_{\delta}$ is a connected curve. 

Then $\B(C; p_1, p_2 , \ldots, p_{\delta})$ is smooth at the point $[(C; p_1, p_2 , \ldots, p_{\delta})]$ of 
dimension $\dim(|L|) - \delta$. 
\end{theorem}
\begin{proof} This follows from \cite[Theorem 3.8]{Tan} since, by our assumptions, the pointed curve 
$(C; p_1 , \ldots, p_{\delta})$ is {\em virtually connected} in the language of \cite[Definition 3.6]{Tan}. 
\end{proof}

For the proof of Theorem \ref{thm:esistenza} we need to recall other fundamental facts.
We first define, for any globally generated
line bundle $L$ of sectional genus $p:=p_a(L) \geq 2$, on a $K3$
surface $S$, and any integer $\delta$ such that $0 < \delta \leq
p-2$, the locus in the Severi variety $V_{|L|, \delta}$,
\begin{equation} \label{eq:viper}
V_{|L|, \delta}^{hyper} :=  \Big \{ C \in V_{|L|, \delta} \; | \;
\mbox{its normalization is hyperelliptic} \Big \}.
\end{equation}
Observe that in particular,
for any $p \geq 3$, one always has $ V_{|L|, p-2}^{hyper} = V_{|L|, p-2} \neq \emptyset$
and, by regularity of $V_{|L|, p-2}$, this is smooth and of dimension two.

Let $\M_g$ be the moduli space of smooth curves of genus $g$, which is quasi-projective
of dimension $3g-3$ for $g \geq 2$. Denote by ${\overline \M}_g$
its Deligne-Mumford compactification. Then ${\overline \M}_g$ is the moduli space of
stable, genus $g$
curves. Let $\Ha_g \subset \M_g$ denote the locus of hyperelliptic curves, which is known to be an
irreducible variety of dimension $2g -1$ (see e.g. \cite {AC}) and
$\overline{\Ha}_g \subset {\overline \M}_g$ be its compactification.

Moreover, recall from \cite[Def.(3.158)]{HM} that a nodal curve $C$ (not necessarily irreducible) is
{\it stably equivalent} to a stable curve $C'$ if $C'$ is obtained from $C$ by contracting to a point all smooth
rational components of $C$ meeting the other components in only one or two points.

As above, we define the degenerated version of $V_{|L|, \delta}^{hyper}$ by
\begin{eqnarray}
\label{eq:defw} W_{|L|, \delta}^{hyper} & := & \Big  \{ C \in W_{|L|, \delta} \; | \;
\mbox{there exists a desingularization $\widetilde{C}$ of $\delta$ of the} \\
\nonumber            & & \mbox{nodes of $C$, such that $\widetilde{C}$ is stably equivalent to a} \\
\nonumber            & & \mbox{(stable) curve $C'$ with $[C'] \in \overline{\Ha}_{p_a(L)-\delta}$} \Big  \}.
\end{eqnarray}

\noindent
Note that, by definition, any such $\widetilde{C}$ is connected.
Similarly as in \eqref{eq:irred}, we have:
\begin{equation} \label{eq:irred2}
  W_{|L|, \delta}^{hyper}  = \cup_{\delta' \geq \delta} V_{|L|, \delta}^{hyper} \;
\mbox{if all the curves in $|L|$ are irreducible}.
\end{equation}

Theorem \ref{thm:esistenza} will be a direct consequence of the next three results,
Propositions \ref{prop:degarg} and \ref{prop:k3speciale} and Lemma \ref{lemma:k3speciale}.
The central degeneration argument is given by the following:

\begin{proposition} \label{prop:degarg}
Let $p \geq 3 $ and $\delta \leq p-2$ be positive integers.
Assume there exists a smooth $K3$ surface $S_0$ with a globally generated, primitive line bundle
$H_0$ on $S_0$ with $p_a(H_0)=p$ and such that ${W}_{|H_0|, \delta}^{hyper}(S_0) \neq \emptyset$ and 
$\dim ({W}_{|H_0|, \delta}^{hyper}(S_0)) \leq 2$.

Then, on the general, primitively marked $K3$ surface $(S,H)$ of genus $p$,
${W}_{|H|, \delta}^{hyper}(S)$ is nonempty and equidimensional of dimension two.
\end{proposition}

\begin{proof}
Let $\B_p$ be the moduli space of primitively marked $K3$ surfaces of genus $p$.
It is well-known that $\B_p$ is smooth and irreducible of dimension $19$, cf. e.g. \cite[Thm.VIII 7.3 and p.~366]{BPV}.
We let $b_0=[(S_0,H_0)] \in \B_p$.
Similarly as in \cite{B1}, consider the {\it scheme of pairs}
\begin{equation}\label{eq:stack}
\mathcal{W}_{p,\delta} := \Big  \{ (S,C) \; | \;  [(S,H)] \in \B_p \; \mbox{and} \;  [C] \in W_{|H|,\delta}(S) \Big  \},
\end{equation}
and the natural projection
\begin{equation} \label{eq:pi}
\pi : \mathcal{W}_{p,\delta} \hpil \B_p.
\end{equation}(The fact that $\mathcal{W}_{p,\delta}$ is a scheme, in fact a locally closed scheme, follows from
the already mentioned proof of  Mumford's theorem on the existence of nodal
rational curves as in
\cite[pp.~351-352]{MM} or \cite[pp.~365-367]{BPV}.)

Note that for general $[(S_b,H_b)]=b \in \B_p$ we have
\[ \pi^{-1}(b)= \cup_{\delta' \geq \delta} V_{|H_b|,\delta'}(S_b) \]
by \eqref{eq:irred} (as $\Pic (S_b) \iso \ZZ[H_b]$), so that $\pi^{-1}(b)$ is nonempty, equidimensional and
of dimension $g:=p-\delta$, by the regularity property recalled above. In particular, $\pi$ is dominant.
Observe that $\mathcal{W}_{p,\delta}$ is singular in codimension one, so in particular it is not normal.

For brevity, let $\mathcal{W}:= \mathcal{W}_{p,\delta}$ and let ${\mathcal C} \stackrel{f}{\to} {\mathcal W} $
be the universal curve.
As in Theorem \ref{thm:sern}, (i) and (ii), in Appendix \ref{A:edoardo},
there exists a commutative diagram
\[
\xymatrix{
{\mathcal C}'  \ar[d]_{f'} \ar[r] & {\mathcal C} \ar[d]^{f} \\
{\mathcal W}_{(\delta)} \ar[r]^{\alpha} &  {\mathcal W},
 }
\]where $\alpha$ is a finite, unramified morphism defining a marking of all the
$\delta$-tuples of nodes of the fibres of $f$
(cf. Theorem \ref{thm:sern}, with $V = {\mathcal W}$, $E_{(\delta)} = {\mathcal W}_{(\delta)}$).
Precisely, by using notation as in Theorem \ref{thm:sern}, if for $w \in {\mathcal W}$ the curve
$ {\mathcal C}(w)$ has $\delta + \tau$ nodes, $\tau \in \ZZ^+$, $\alpha^{-1}(w)$
consists of $ \delta + \tau \choose \delta$ elements, since any $\eta_w \in \alpha^{-1}(w)$ parametrizes
an unordered, marked $\delta$-tuple of the $\delta + \tau$ nodes of ${\mathcal C}(w)$.

Let $\eta_w \in {\mathcal W}_{(\delta)}$. Then $\eta_w$ is represented by a pointed curve $(C; p_1, p_2 , \ldots, p_{\delta})$, where $(S,C) \in  {\mathcal W}$ and where $p_1, p_2 , \ldots, p_{\delta}$ are 
$\delta$ marked nodes on $C$. 

Let ${\mathcal W}(S,H)$ (resp. ${\mathcal W}_{(\delta)}(S,H)$) be the fibre of $\pi$ (resp. of $\alpha \circ \pi$) 
over $[(S,H)] \in \B_p$, and let 
$$
\alpha(S,H) : {\mathcal W}_{(\delta)}(S,H) \hpil {\mathcal W}(S,H)
$$
be the induced morphism. For $\eta_w \in {\mathcal W}_{(\delta)}(S,H)$ as above, we have
\begin{equation}\label{eq:edo}
T_{[\eta_w]}({\mathcal W}_{(\delta)}(S,H) ) \cong T_{[(C; p_1, p_2 , \ldots, p_{\delta})]}(\B(C; p_1, p_2 , \ldots, p_{\delta}) ),  
\end{equation}where $\B(C; p_1, p_2 , \ldots, p_{\delta})$ is as in \eqref{eq:branch}. Indeed, since $\alpha$ is finite and unramified, then also $\alpha(S,H)$ is. Therefore, it suffices to consider the image of the differential $d \alpha(S,H)_{[\eta_w]}$. The latter is given by first-order deformations of $C$ in $S$ (equivalently in $|H|$) 
that are locally trivial at the $\delta$ marked nodes; these are precisely given by $T_{[(C; p_1, p_2 , \ldots, p_{\delta})]}(\B(C; p_1, p_2 , \ldots, p_{\delta}) )$ (cf. \cite[Remark 3.5]{Tan}).

Let $\widetilde{\mathcal{W}}_{(\delta)}$ be the smooth locus of $\mathcal{W}_{(\delta)}$. By Theorem \ref{thm:Tan3} 
and by \eqref{eq:edo}, together with the fact that $\B_p$ is smooth, 
$\widetilde{\mathcal{W}}_{(\delta)}$ contains all the pairs $(S,C)$ with $\delta$ marked nodes on $C$, such that 
$|C|$ is globally generated (i.e. its general element is a smooth, irreducible curve) and the partial normalization 
of $C$ at these marked nodes is a connected curve. More precisely,  by the proof of Mumford's theorem on the existence of nodal rational curves on $K3$ surfaces, as in \cite[pp.~351-352]{MM} or \cite[pp.~365-367]{BPV}), any irreducible component of $\mathcal{W}_{(\delta)}$ has dimension $\geq 19+p-\delta=19+g$; furthermore, by \eqref{eq:edo}, $\dim(T_{[\eta_w]}(\mathcal{W}_{(\delta)}(S,H))) = g$, where $\eta_w$ represents $(S,C)$ with $C$ with the $\delta$ 
marked nodes. It also follows that $\mathcal{W}_{(\delta)}$ is smooth, of dimension $19+g$ at these points.

If we restrict ${\mathcal C}'$ to $\widetilde{\mathcal{W}}_{(\delta)}$, from Theorem \ref{thm:sern}, (iv) and (v),
we have a commutative diagram
\[
\xymatrix{
\widetilde{\mathcal C}  \ar[d]_{\widetilde{f}} \ar[r] & {\mathcal C} \ar[d]^{f} \\
\widetilde{\mathcal{W}}_{(\delta)} \ar[r]^{\widetilde{\alpha}} &  {\mathcal W},
 }
\]where $\widetilde{\alpha} = \alpha|_{\widetilde{\mathcal{W}}_{(\delta)}}$ and where
$\widetilde{f}$ is the flat family of partial normalizations at $\delta$ nodes of the curves
parametrized by $\alpha(\widetilde{\mathcal{W}}_{(\delta)})$ (in the notation of Theorem \ref{thm:sern} in Appendix 
\ref{A:edoardo}, $\widetilde{f} = \overline{f}$ in (v) and $\widetilde{\mathcal C} = \overline{\mathcal C}$ in 
(iii) and (iv)).

There is an obvious rational map
\[
\xymatrix{ \widetilde{\mathcal{W}}_{(\delta)}
\ar@{-->}[r]^{c} &  \overline{\M}_g, 
}
\]
defined on the open dense subscheme
$\widetilde{\mathcal{W}}_{(\delta)}^0 \subset \widetilde{\mathcal{W}}_{(\delta)}$ 
such that, for $\eta_w \in
\widetilde{\mathcal{W}}_{(\delta)}^0 $,
$\widetilde{\mathcal C}(\eta_w)$ is  stably equivalent to a stable curve of genus $g$.

Set $\psi := c|_{\widetilde{\mathcal{W}}_{(\delta)}^0}$. By definition, for any
$\eta_w \in \widetilde{\mathcal{W}}_{(\delta)}^0$,
the  map $\psi$ contracts all possible smooth rational components of $\widetilde{\mathcal C}(\eta_w)$
meeting the other components in only one or two points and maps 
the resulting stable curve into its equivalence class in $\overline{\M}_g$.

Pick any $C_0 \in W_{|H_0|, \delta}^{hyper}(S_0)$ and let $w_0 = [(S_0, C_0)] \in \mathcal{W}$ be the corresponding point. 
Now $|H_0|$ is globally generated and the normalization of $C_0$ at some $\delta$ nodes satisfying the conditions 
in \eqref{eq:defw} is a connected curve. Therefore, letting $\eta_{w_0} \in \alpha^{-1}(w_0)$ be the point 
corresponding to marking these $\delta$ nodes, we have that 
$\eta_{w_0} \in \widetilde{\mathcal{W}}_{(\delta)}^0$ and the map $c$ is defined at $\eta_{w_0}$.

Let $\widetilde{\V} \subseteq \widetilde{\mathcal{W}}_{(\delta)}^0$ be the irreducible component containing 
$\eta_{w_0}$; then, as proved above, $\dim ( \widetilde{\V} ) = 19 + g$. 

By assumption, $\psi(\widetilde{\V}) \cap \overline{\Ha}_g \neq \emptyset$.
Hence, for any irreducible component
${\mathcal K} \sub \psi(\widetilde{\V}) \cap \overline{\Ha}_g$, we have
\begin{equation} \label{eq:main(*)}
\dim ({\mathcal K})  \geq \dim (\psi(\widetilde{\V})) +
\dim (\overline{\Ha}_g) - \dim (\overline{\M}_g) = \dim (\psi(\widetilde{\V})) + 2-g.
\end{equation}

Pick any ${\mathcal K}$ containing $\psi(\eta_{w_0})$ and let
$\Ii \sub \psi_{|{\widetilde{\V}}}^{-1}({\mathcal K})$ be any irreducible component containing
$\eta_{w_0}$. Since the general fibre of $\psi_{|_{\widetilde{\V}}}$ has dimension
$\dim (\widetilde{\V}) - \dim (\psi(\widetilde{\V})) = 19+g-\dim (\psi(\widetilde{\V}))$, 
from \eqref{eq:main(*)} we have
\begin{eqnarray} \label{eq:piu21}
 \dim (\Ii)  & = & \dim ({\mathcal K}) + 19+g-\dim (\psi(\widetilde{\V})) \\
\nonumber & \geq & \dim (\psi (\widetilde{\V})) + 2-g + 19 + g - \dim (\psi(\widetilde{\V})) = 21.
\end{eqnarray}

Consider now
\begin{equation} \label{eq:base}
\pi \circ (\widetilde{\alpha}|_{\Ii}) : \Ii \longrightarrow \B_p.
\end{equation}
Since, by assumption, the fibre over $b_0 = [(S_0, H_0)]$ is at most two-dimensional, we conclude from
\eqref{eq:piu21} that $\pi \circ (\widetilde{\alpha}|_{\Ii}) $ is dominant, that
all the fibres are precisely two-dimensional and that $\dim (\Ii)=21$.
This shows that ${W}_{|H|, \delta}^{hyper} \neq \emptyset$ for
general $[(S,H)] \in \B_p$ and
Lemma \ref{cor:expdim} implies that in fact any irreducible
component of ${W}_{|H|, \delta}^{hyper}(S)$ has 
dimension two.
\end{proof}

\begin{remark}\label{rem:p=4,5}{\normalfont
In particular, Lemma \ref{cor:expdim}, Proposition \ref{prop:degarg} 
and \cite[Examples 2.8 and 2.10]{fkp} prove Theorem \ref{thm:esistenza} for $p=4$ and $5$.
}
\end{remark}

We next construct the desired special primitively marked $K3$ surface:

\begin{prop} \label{prop:k3speciale}
Let $d \geq 2$ and $k \geq 1$ be integers. There exists a $K3$ surface $S_0$ with
\[ \Pic (S_0) = \ZZ [E] \+ \ZZ [F] \+ \ZZ [R]  \]
and intersection matrix
\[  \left[
  \begin{array}{cccc}
   E^2        &   E.F      &  E.R               \\
   F.E        &   F^2      &  F.R                   \\
   R.E        &   R.F      &  R^2
    \end{array} \right]  =
    \left[
  \begin{array}{ccccc}
  0 & d & k       \\
  d & 0 & k         \\
  k & k & -2
    \end{array} \right],     \]
and such that the following conditions are satisfied:
\begin{itemize}
\item[(a)] $|E|$ and $|F|$ are elliptic pencils;
\item[(b)] $R$ is a smooth, irreducible rational curve.
\item[(c)] $H_0:=E+F+R$ is globally generated, in particular the general member of $|H_0|$ is a smooth, irreducible
curve of arithmetic genus $p:=2k+d$;
\item[(d)] the only effective decompositions of $H_0$ are
\[ H_0 \sim E + F + R \sim (E+F) + R \sim (E+R) + F \sim (F+R)+E.\]
\end{itemize}
\end{prop}

\begin{proof} Since the lattice  has signature $(1,2)$, then, by a result of Nikulin \cite{nik}
(see also \cite[Cor. 2.9(i)]{Morr}), there is a
$K3$
surface $S_0$ with that as Picard lattice.
Performing
Picard-Lefschetz reflections on the lattice, we can assume that $H_0$ is nef,
by \cite[VIII, Prop. 3.9]{BPV}. Straightforward calculations on the Picard
lattice rules out the existence of effective divisors $\Gamma$ satisfying
$\Gamma^2=-2$ and $\Gamma.E <0$ or $\Gamma.F <0$, or $\Gamma^2=0$ and
$\Gamma.H_0=1$. Hence (a) and (c) follow from \cite[Prop. 2.6 and (2.7)]{SD}. Similarly one computes that if
$\Gamma >0$, $\Gamma^2=-2$ and $\Gamma.R <0$, then $\Gamma=R$, proving (b).

Similarly, (d) is proved by direct calculations using the nefness of $E$, $F$ and $H_0$ and recalling that by Riemann-Roch and Serre duality a divisor $D$ on a $K3$ surface is effective and irreducible only if $D^2 \geq -2$ and $D.N >0$ for some nef divisor $N$.
\end{proof}

The following result, together with \eqref{eq:irred2} and
Proposition \ref{prop:degarg}, now concludes the proof of Theorem
\ref{thm:esistenza}  and Corollary \ref{cor2:g=3}. From Remark
\ref{rem:p=4,5}, we need only consider $p \geq 6$.

\begin{lemma} \label{lemma:k3speciale}
Let $p \geq 6 $ be an  integer. There exists a smooth $K3$ surface $S_0$ with a globally generated, primitive line bundle
$H_0$ on $S_0$ with $p=p_a(H_0)$ such that
\begin{itemize}
\item[(a)] $W_{|H_0|, p-3}^{hyper}(S_0) \neq \emptyset$;
\item[(b)] $\dim (W_{|H_0|, p-3}^{hyper}(S_0)) = 2$;
\item[(c)] there exists a component of $W_{|H_0|, p-3}^{hyper}(S_0)$ whose general member
deforms to a curve $[C_t] \in
V_{|H_t|, p-3}^{hyper}(S_t)$, for general $[(S_t,H_t)] \in \B_p$;
\item[(d)] for general $[(S_t,H_t)] \in \B_p$, the two-dimensional irreducible component
$V_t \sub V_{|H_t|, p-3}^{hyper}(S_t)$ given by (c), satisfies $\dim (R_{V_t})=3$
(with notation as in \S\;\ref{SS:BB}).
\end{itemize}
\end{lemma}

\begin{proof}
Set $k=1$ if $p$ is even and $k=2$ if $p$ is odd and let $d:=p-2k \geq 2$.
Consider the marked $K3$ surface
$(S_0, H_0)$ in Proposition \ref{prop:k3speciale}.

We will consider two general smooth elliptic curves $E_0 \in |E|$ and
$F_0 \in |F|$ and curves of the form
\[
C_0 := E_0 \cup F_0 \cup R, 
\]with transversal intersections and a desingularization 
\begin{equation}\label{eq:czero}
\widetilde{C}_0 = \widetilde{E}_0 \cup \widetilde{F}_0 \cup \widetilde{R} \khpil C_0 
\end{equation}of the $\delta:=p-3=d+2k-3$ nodes marked in Figure \ref{fig:31} below, that is, all but one of each
of the intersection points $E_0 \cap F_0$,
$E_0 \cap R$ and $F_0 \cap R$.

\begin{figure}[ht]
\[
\includegraphics{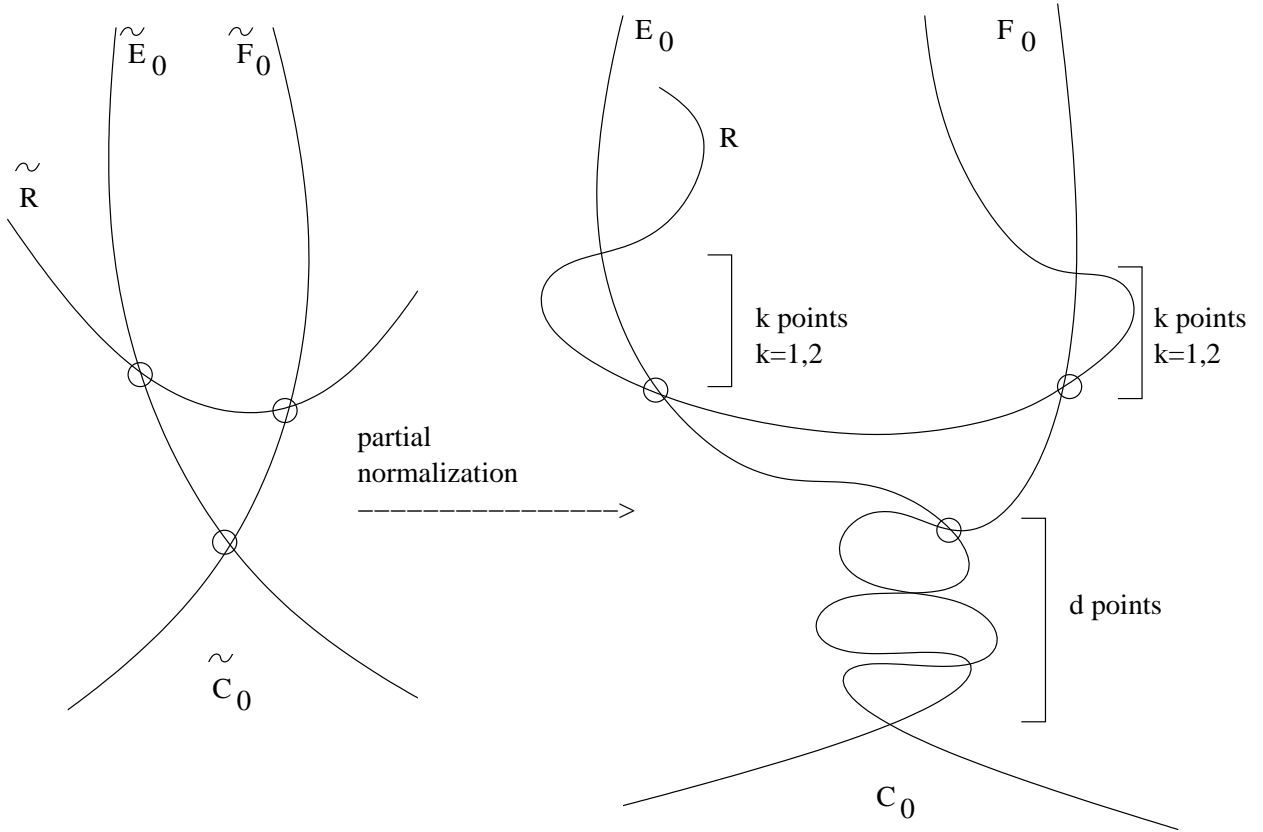}
\]
\caption{The curves $C_0$ and $\widetilde{C}_0$}
\label{fig:31}
\end{figure}

\noindent Then
$[C_0] \in W_{|H_0|,p-3}^{hyper}$, as $\widetilde{C}_0$ is stably equivalent to a union of two smooth
elliptic curves intersecting in two points (cf. \cite[Exercise (3.162)]{HM}),
proving (a).
Clearly the closure of the family we have constructed is isomorphic
to $|E| \x |F| \iso \PP^1 \x \PP^1$, and is therefore two-dimensional.
Denote by $W_0 \subset W_{|H_0|, p-3}^{hyper}$ this two-dimensional subscheme.

We will now show that any irreducible component $W$ of $W_{|H_0|,p-3}^{hyper}$
has dimension $\leq 2$.

A central observation, which will be used together with Theorem \ref{thm:apppg}(b), will be that, with the above choices of $k$, we have
\begin{equation} \label{eq:parita}
 E.H_0=F.H_0=d+k=p-k \; \mbox{is odd}.
\end{equation}

We start by considering families of reducible curves. These are all classified in
Proposition \ref{prop:k3speciale}(d).

If the general element in $W$ is of the form $D \cup R$, for $D
\in |E+F|$, then in order to have a partial desingularization
$\widetilde{D} \cup \widetilde{R}$ to be (degenerated)
hyperelliptic, we must have $\deg (\widetilde{D} \cap
\widetilde{R})=2$, so that we must desingularize $2(k-1)$ of the
intersection points of $D \cap R$. Finally, as $p_a(\widetilde{D}
\cup \widetilde{R})=3$, we must have $p_a(\widetilde{D})=2$.
Therefore $W \sub W_D \x \{R\} \iso W_D$, where $W_D \subset|D|$ is a
subfamily of irreducible curves of geometric genus $\leq 2$. It
follows that $\dim (W) \leq \dim (W_D) \leq 2$, by Lemma \ref{lemma:dimre}.

If the general element in $W$ is of the form $D \cup E$, for $D \in |F+R|$, then in order to have
a partial desingularization $\widetilde{D} \cup \widetilde{R}$ that is (degenerated) hyperelliptic,
we must have $\deg (\widetilde{D} \cap \widetilde{E})=2$. If the projection $W \khpil |E|$ is dominant,
this means that $\mathfrak{g}^1_2(\widetilde{D}) \sub |f^*E|_{|\widetilde{D}}$, where
$f:\widetilde{S} \khpil S$ denotes the composition of blow-ups of $S$ that induces
the partial desingularization $\widetilde{D} \cup \widetilde{R} \khpil D \cup R$.
But this would mean that $|f^*E|_{|\widetilde{D}}$, which is base point free
on $\widetilde{D}$, is composed with the
$\mathfrak{g}^1_2(\widetilde{D})$, a contradiction, as $\deg (\cO_{\widetilde{D}}(f^*E))=E.D=E.H_0$
is odd by \eqref{eq:parita}.
Therefore, the projection $W \khpil |E|$ is not dominant, whence
$\dim (W) \leq \dim (|D|) = \frac{1}{2}D^2+1= k \leq 2$, as
desired.
By symmetry, the case where the general element in $W$ is of the form $D \cup F$, for $D \in |E+R|$
is treated in the same way.

Finally, we have to consider the case of a family $W \sub |H_0|$ of irreducible curves.

In this case assume $\dim (W) \geq
3$, and let $C$ be a general curve parametrized by $W$. Then by Theorem \ref{thm:apppg} (b), there exists an
 effective decomposition
into moving classes $H_0 \sim M+N$ such that
\[ \mathfrak{g}^1_2(\widetilde{C}) \sub |f^*M|_{|\widetilde{C}}, \; |f^*N|_{|\widetilde{C}},\]
where $f:\widetilde{S} \khpil S$ denotes the succession of blow ups of $S$ that induces the normalization $\widetilde C\to C$.
From Proposition \ref{prop:k3speciale}(d)
we see that we must have
\[ \mathfrak{g}^1_2(\widetilde{C}) \sub |f^*E|_{|\widetilde{C}}, \; \mbox{or} \; |f^*F|_{|\widetilde{C}},\]
which means that either $|f^*E|_{|\widetilde{C}}$ or
$|f^*F|_{|\widetilde{C}}$
 is composed with the $\mathfrak{g}^1_2(\widetilde{C})$, again a contradiction,
as both have odd degree by \eqref{eq:parita}. We have therefore proved (b).

To prove (c) we will show that any $[C_0] \in W_{|H_0|,p-3}^{hyper}$ in the two-dimensional, irreducible component
$W_0$ considered above in fact deforms to a curve $[C_t] \in
W_{|H_t|, p-3}^{hyper}(S_t)$, for general $[(S_t,H_t)] \in \B_p$, that has precisely
$\delta=p-3$ nodes (cf. \eqref{eq:irred2}).

To this end, denote by $\mathcal{S} \khpil \B_p$ the universal family of $K3$ surfaces, 
$\tilde{f}: \widetilde{\C} \khpil \widetilde{\W}_{(\delta)}$ and
$\I \subset \widetilde{\W}_{(\delta)}$ as in the proof of Proposition \ref{prop:degarg}, and let
$\varphi: \widetilde{\C}_{\Ii} \khpil \Ii$ be the restriction of $\tilde{f}$.

Since the fiber over $[(S_0, H_0)]$ of ${\Ii} \khpil \B_p$ as in \eqref{eq:base}
contains an open, dense subset of $\PP^1 \x \PP^1$, we can find a smooth, irreducible curve $B \subset \Ii$ satisfying: for $x \in B$ general,
$\varphi^{-1}(x)$ is a (partial) desingularization of $\delta=p-3$ of the nodes of a curve in
$W_{|H_t|,\delta}(S_t)$ (cf. \eqref{eq:w}), for general $[(S_t, H_t)] \in \B_p$, and 
$\varphi^{-1}(x) \in \overline{\mathcal{H}}_3 \subset \overline{\M}_3$; moreover $B$ contains a point
$x_0 \in \Ii$ such that $\varphi^{-1}(x_0)$ is $\widetilde{C}_0$ as in \eqref{eq:czero}, for $C_0$ general
in $W_0$.

Let $\varphi_B : \widetilde{{\mathcal C}}_B \to B$ be the 
induced universal curve. Since the dualizing sheaf of $\varphi_B^{-1} (x_0)=\widetilde{C}_0$ is globally generated (as each component intersects the others in two points), we in fact have, possibly after substituting $B$ with an open neighbourhood of $x_0$,
a morphism 
$\gamma_B: \widetilde{\C}_B \khpil \PP(\tilde{\varphi}_* (\omega_{\widetilde{\C}/B}))$
over $B$ that is $2:1$ on the general fibre ${\varphi}_B^{-1}(x)$ and contracts the 
rational component
$\widetilde{R}$ of  $\varphi_B^{-1}(x_0)$ and maps the two elliptic curves
$\widetilde{E}_0$ and $\widetilde{F}_0$ each $2:1$ onto (different) $\Pp^1$s (cf. \eqref{eq:czero}
and Figure \ref{fig:31}).

Let $\nu: \widetilde{{\mathcal C}}'_B \khpil \widetilde{\C}_B$ be the normalization and
\[
\xymatrix{ \widetilde{\C}'_B \ar[r]^{\gamma_1} &  \widetilde{\C}''_B \ar[r]^{\hspace{-1cm}\gamma_2} &
\PP(\tilde{\varphi}_* (\omega_{\widetilde{\C}_B/B}))
}
\]
the Stein factorization of $\gamma_B \circ \nu$. In particular,
$\gamma_2$ is finite of degree two onto its image.  Moreover, $\nu \circ \varphi_B: \widetilde{\C}'_B \to B$ is a flat family whose general fiber $(\nu \circ \varphi_B)^{-1}(x)$ is a desingularization of 
$\varphi_B^{-1}(x) \in \widetilde{\C}_B$. Let $p_g$ be the geometric genus of this general fibre.

Let ${\mathcal D} \subset \widetilde{\C}'_B$ be the strict transform via $\gamma_1$ of the
closure of the branch divisor of $\gamma_2$ on the smooth locus of
$\widetilde{\C}''_B$.
By Riemann-Hurwitz, for general $x \in B$, we have
${\mathcal D}.\varphi_B^{-1}(x)=2p_g+2$, whereas
${\mathcal D}.\varphi_B^{-1} (x_0) \geq 8$, as
the curve $\gamma_1(\varphi_B^{-1}(x_0))$ contains two smooth elliptic curves,
each being mapped $2:1$ by $\gamma_2$ onto (different) $\PP^1$s.
This implies $p_g=3$.  Since, for general $x \in B$, we have $p_g \leq p_a(\varphi_B^{-1}(x))=p-\delta=3$, we find that $\varphi_B^{-1}(x)$ is smooth. This means that the general curve in $W_{|H_t|,\delta}(S_t)$, for $(S_t,H_t) \in \B_p$ general, has precisely $\delta=p-3$ nodes, 
proving (c).

To prove (d), again we consider the morphism (up to possibly restricting $\I$ as above) 
\[ \gamma_{\Ii}: \C_{\Ii} \to \PP(\varphi_* (\omega_{\C_{\Ii}/\Ii}))
\]over $\Ii$ which, apart some possible contractions of rational components in special fibres
over $\Ii$, is relatively $2:1$ onto its image.
We have a natural morphism $h : \C_{\Ii} \to \mathcal{S}$, inducing a natural map
\[
\Phi: \im (\gamma_{\Ii}) - - \to \Sym^2(\mathcal{S}), \]
whose domain has nonempty intersection with every fibre over $\B_p$.

Let $\mathcal{R}:= \overline{\im (\Phi)}$. Then $\mathcal{R} \cap
\Sym^2(S_t)=R_{V_t}$, for general $[(S_t,H_t)] \in \B_p$.
One easily sees that
\[ \{\Sym^2(E')\}_{E' \in |E|} \cup \{\Sym^2(F')\}_{F' \in |F|} \subseteq
\overline{\mathcal{R} \cap \Sym^2(S_0)}. \]
Since the two varieties on the left are threefolds, we have $\dim (\Phi^{-1} (\xi_0)) =0$
for general $\xi_0 \in \mathcal{R} \cap \Sym^2(S_0) \subset \mathcal{R}$.
Therefore, for general $\xi \in \mathcal{R}$, we have $\dim
(\Phi^{-1} (\xi)) =0$, so that $\dim
(\mathcal{R})=\dim(\C_{\Ii}) = \dim
(\Ii)+1=22$, whence $\dim (R_{V_t})=22- \dim
(\B_p)=3$.
\end{proof}

\begin{remark} \label{rem:g=g0}
{\rm For general $[(S_t, H_t)] \in \B_p$ the obtained curves in the last proof have in fact $\delta=p-3$}
non-neutral {\rm nodes (cf. \cite[\S 3]{fkp}). In fact a desingularization of less than $p-3$ nodes of $C_t$ admits no $\mathfrak{g}^1_2$s, as clearly a desingularization of less than $p-3$ nodes of $C_0$ is not stably equivalent to a curve in the hyperelliptic locus $\overline{\Ha}_3 \subset \overline{\M}_3$.}
\end{remark}

%
\section{On the Mori cone of the Hilbert square of a $K3$ surface} \label{S:K3}
%

In this section we first summarize central results on the Hilbert square of a $K3$ surface and show how to compute the class of a rational curve in $\ss$. Then we discuss the relations between the existence of curves on $S$ and the slope of the Mori cone of $\ss$, that is, the cone of effective classes in 
$N_1(\ss)_{\RR}$. In particular, we show how to deduce the bound \eqref{eq:g=3boundinf} from
Theorem \ref{thm:esistenza} and \eqref{eq:seshadribound} from known results about Seshadri constants.
Finally, we discuss 
the relation between the
existence of a curve on $S$ with given singular Brill-Noether number 
and the slope of the Mori cone of $\ss$.

%
\subsection {Preliminaries on $\ss$ for a $K3$ surface} \label{S:prel}
%

Recall that for any smooth surface $S$ we have
\begin{equation} \label{eq:decH2}
H^2(S^{[2]}, \mathbb Z) \iso H^2(S, \mathbb Z) \oplus \mathbb Z \mathfrak{e},
\end{equation}
where $\Delta:=2\mathfrak{e}$ is the class of the divisor parametrizing
$0$-dimensional subschemes supported on a single point (see
\cite{B}). So we may identify a class in $H^2(S, \mathbb Z)$ with
its image in $H^2(S^{[2]}, \mathbb Z)$. When $S$ is a $K3$ surface the cohomology group
$H^2(S^{[2]}, \mathbb Z)$ is endowed with a quadratic form $q$,
called the {\em Beauville-Bogomolov form}, such that its
restriction to $H^2(S, \mathbb Z)$ is simply the cup product on
$S$, the two factors $H^2(S, \mathbb Z)$ and $\mathbb Z \mathfrak{e}$ are
orthogonal with respect to this form and $q(\mathfrak{e})=-2$. The decomposition
\eqref{eq:decH2} induces an isomorphism
\begin{equation} \label{eq:decpic}
\Pic (S^{[2]}) \iso \Pic (S)  \oplus \mathbb Z [\mathfrak{e}],
\end{equation}
and each divisor $D$ on $S$ corresponds to the divisor on $S^{[2]}$, by abuse of notation also denoted by $D$,
consisting of length-two subschemes with some support on $D$.

Given a primitive class $\alpha\in H_2(S^{[2]}, \mathbb Z)$, there exists a unique class $w_{\alpha}\in
H^2(S^{[2]}, \mathbb Q)$ such that $\alpha.v=q(w_{\alpha},v)$, for all $v\in
H^2(S^{[2]}, \mathbb Z)$, and one sets
\begin{equation}\label{eq:qalfa}
q(\alpha):= q(w_{\alpha}).
\end{equation}We denote also by $\rho_{\alpha}\in H^2(S^{[2]}, \mathbb Z)$ the
corresponding primitive $(1,1)$-class such that $\rho_{\alpha}=cw_{\alpha}$, for some $c>0$ (for further details, we refer the reader to \cite{HT}). 

If now $\Pic (S) = \ZZ[H]$,
then the N\'eron-Severi group of $\ss$ has rank two. We may take as
generators
of $N_1(\ss)_{\RR}$ the class $\PP^1_{\Delta}$ of a rational curve in the ruling of
the exceptional
divisor $\Delta \subset \ss$, and the class of the curve in $\ss$ defined as
follows
\[
 \{\xi\in \ss | \Supp(\xi)=\{p_0, y\} \; | \;  y\in Y\},
\]
where $Y$ is a curve in $|H|$ and $p_0$ is a fixed point on $S$. By abuse
of notation, we still denote the class of the curve in $\ss$ by $Y$.
 Note that we always have that
\begin{equation} \label{P1sulbordo}
 \PP^1_{\Delta} \; \mbox{lies on the boundary of the Mori cone.}
\end{equation}
Indeed, the curve $\PP^1_{\Delta}$ is contracted by the Hilbert-Chow morphism
$\ss\to \Sym^2(S)$, so that the pull-back of an ample divisor on
$\Sym^2(S)$
is nef, but zero along $\PP^1_{\Delta}$.

Therefore,
describing
the Mori cone $\NE (\ss)$ amounts, by (\ref{P1sulbordo}), to computing
\begin{equation}\label{inf}
 \slope(\NE (\ss)):= \inf \Big \{\frac {a}{b} \; | \:  aY - b\PP^1_{\Delta}\in N_1(\ss)\ \textrm{is
effective,}\; a,b \in \QQ^+ \Big \}.
\end{equation}
We will also call the (possibly infinite) number $a/b$ associated to an irreducible curve 
$X \sim_{alg} aY-b\PP^1_{\Delta}$ with $a >0$ and $b \geq 0$, the {\it slope of the curve} $X$ and denote it by $\slope(X)$. Thus, 
the smaller 
$\slope(X)$ is, the nearer is $X$ to the boundary of $\NE (\ss)$. 

By a general result due to Huybrechts \cite[Prop. 3.2]{H2} and Boucksom \cite{Bou}, a divisor $D$ on $\ss$ is ample
if and only if $q(D)>0$ and $D.R>0$ for any (possibly
singular) rational curve
$R\subset \ss$. As a consequence, if the Mori cone is closed then the boundary (which remains to be determined)  
is generated by the class of a rational curve (the other boundary is generated by 
$\PP^1_{\Delta}$, by \eqref{P1sulbordo}). This means that one would have $\slope(\NE (\ss))=\slope_{rat}(\NE (\ss))$, 
where 
\begin{equation} \label{inf-rat}
 \slope_{rat}(\NE (\ss)):= \inf \Big \{\frac {a}{b} \; | \;  aY - b\PP^1_{\Delta}\in N_1(\ss) \; \mbox{is
the class of a rational curve,} \; a,b \in \QQ^+ \Big \}.
\end{equation}
(A priori, one only has
$\slope(\NE (\ss)) \leq \slope_{rat}(\NE (\ss))$.)

Hassett and Tschinkel \cite{HT} make
a precise prediction on the geometric and numerical properties of such extremal rational
curves in $\ss$. Indeed, according to their conjectures \cite[p.~1206 and Conj. 3.6]{HT}, the extremal ray $R$ has to be generated either by the class of a line inside a $\PP^2$, such that $q(R) = - \frac{5}{2}$ as in \eqref{eq:qalfa}, or by the class of a rational curve that is a fibre of a $\PP^1$-bundle over a $K3$ surface and 
such that $q(R)=-2$ or $-\frac{1}{2}$.

%
\subsection{The classes of rational curves in $\ss$} \label{S:class}
%

Assume that $\Pic (S) = \ZZ[H]$ with $p_a(H)=p_a \geq 2$. Let $X \subset \ss$ be an irreducible rational curve. Let $C_X \subset S$ be the corresponding curve as in \S\;\ref{S:corr} and assume that
$C_X \in |mH|$ with $m \geq 1$. (In particular, $m \geq 2$ if we are in case (II)).
We can write
\[
 X \sim_{alg} a_1 Y + a_2 \PP^1_{\Delta}.
\]
Since $X.H=m(2p_a-2)$, $Y.H = 2p_a-2$ and $\PP^1_{\Delta}.H=0$
 by the very definition of $H$ as a divisor in $\ss$, and
 $Y.\mathfrak{e}=0$ and $\PP^1_{\Delta}.\mathfrak{e}=-2$, we obtain, defining
$g_0(X):= X.\mathfrak{e}-1$, 
\begin{equation}\label{eq:classX}
 X \sim_{alg} mY - \Big( \frac{g_0(X)+1}{2} \Big ) \PP^1_{\Delta}.
\end{equation}

To compute $g_0(X)$, consider the diagram \eqref{eq:diarami}. Since 
$\nu_X^* \cO_{X}(\Delta) \iso (\nu_X^* \cO_{X}(\mathfrak{e}))^{\* 2}$, the double cover
$f$ is defined by $\nu_X^* \cO_{X}(\Delta)$. By Riemann-Hurwitz we therefore get
\begin{equation} \label{eq:g0}
g_0(X)=p_a(\widetilde{C}_X).
\end{equation}
Note that in the cases (II) and (III) in the correspondence in \S\;\ref{S:corr},
$X.\mathfrak{e}=g_0(X)+1$ is precisely the length of the intersection scheme 
$\widetilde{C}_{X,1} \cap \widetilde{C}_{X,2}$, where $\widetilde{C}_{X}= 
\widetilde{C}_{X,1} \cup \widetilde{C}_{X,2}$. In case (III), since
$\tilde{\nu}: \widetilde{C}_{X} \to S$ contracts one of the two components of
$\widetilde{C}_{X}$ to a point $x_X \in S$, we obtain that 
\begin{equation} \label{eq:g0'}
g_0(X)=\mult_{x_X}(C_X)-1 \; \mbox{(if $C_X$ is of type (III))}.
\end{equation}

One can check that for all divisors $D$ in $\ss$ one has
$X.D = q(w_X, D)$ with 
\begin{equation}\label{eq:w_X}
 w_X:= mH - \Big ( \frac{g_0(X)+1}{2} \Big ) \mathfrak{e}  \in H^2(\ss, \QQ).
\end{equation}
In particular, $2w_X \in H^2(S^{[2]}, \mathbb Z)$.

From \eqref{inf} and \eqref{eq:classX}  we see that
searching for irreducible rational curves in (or at least ``near'') the boundary of the Mori cone of $\ss$, or with negative square $q(X)$, amounts to searching for irreducible
curves in $|mH|$ with (partial) hyperelliptic normalizations of high genus (case (I)), or to irreducible rational curves in $|mH|$ with high multiplicity at a point (case (III)), or to irreducible rational curves on $S$ with some correspondence between some coverings of their normalizations (case (II)).
Moreover, we should search for curves with as low $m$ as possible. Now $m \geq 2$ in case (II), as remarked above. Moreover, any rational curve in $|H|$ on a general $S$ is nodal, by a result of Chen \cite[Thm. 1.1]{Ch2}
(the same is also conjectured for rational curves in $|mH|$ for $m >1$, see \cite[Conj. 1.2]{Ch1}),
so that $g_0(X) \leq 1$ if $C_X$ is of type (III) in these cases, by \eqref{eq:g0'}. Hence,
we see that the most natural candidates are irreducible curves in $|H|$ with hyperelliptic normalizations. 

By the above, an irreducible curve $C \in |mH|$ with hyperelliptic normalization defines, by the unicity of the 
$\g^1_2$, a unique irreducible 
rational curve $X=R_C \subset \ss$ with class 
\begin{equation}\label{eq:classR}
 R_C \sim_{alg} mY - \Big( \frac{g_0(C)+1}{2} \Big ) \PP^1_{\Delta},
\end{equation}
where $g_0(C):=g_0(R_C)$ is well-defined as
\begin{equation} \label{eq:defg0}
g_0(C):=\mbox{the arithmetic genus of a minimal partial desingularization of $C$ admitting a $\g^1_2$.}
\end{equation}
(For example, if $C$ is nodal,
then we simply take the desingularization of the  
{\it non-neutral nodes} of $C$, cf. \cite[\S 3]{fkp}). 
From \eqref{inf} we then get
\begin{equation} \label{eq:slope}
 \slope(\NE (\ss)) \leq \frac{2m}{g_0(C)+1} \leq \frac{2m}{p_g(C)+1}, \; \mbox{if there exists a $C \in |mH|$ with hyp. norm.}
\end{equation}
and, by \eqref{eq:qalfa} and \eqref{eq:w_X},
\begin{equation}\label{eq:square}
 q(R_C)= 2m^2(p_a-1) - \frac{(g_0(C)+1)^2}{2} \leq 2m^2(p_a-1) - \frac{(p_g(C)+1)^2}{2} .
\end{equation}
In particular, the higher $g_0(C)$ (or $p_g(C)$) is - thus the more ``unexpected'' the curve on $S$ is from a Brill-Noether theory point of view - the lower is the bound on the slope of $\NE(\ss)$ and the more negative is the square $q(R_C)$ in $\ss$.

%
\subsection{The invariant $\rho_{sing}$, Seshadri constants, the ``hyperelliptic existence 
problem'' and the slope of the Mori cone} \label{S:HEP}
%

In \cite{fkp} we introduced a singular Brill-Noether  invariant
\begin{equation} \label{eq:singBN}
\rho_{sing}(p_a,r, d, g):=\rho(g, r, d) +p_a- g,
\end{equation}
in order to study linear series on the normalization of singular
curves.  Precisely, we proved 

\begin{thm} \label{thm:ijm}
Let $S$ be a $K3$ surface such that $\Pic (S) \iso \ZZ[H]$ with $p_a:=p_a(H) \geq 2$. Let $C \in |H|$ 
and $\widetilde{C} \to C$ be a  partial normalization of $C$, such that $g:=p_a(\widetilde{C})$.

If $\rho_{sing}(p_a,r, d, g) <0$, then $\widetilde{C}$ carries no $\g^r_d$.
\end{thm}

\begin{proof}
One easily sees that the proof of \cite[Thm. 1]{fkp} also holds for a {\it partial} normalization
of $C$.
\end{proof}

 For $r=1$ and $d=2$, 
we have
\begin{equation}\label{rho<0}
\rho_{sing}(p_a,1, 2, g)<0 \Leftrightarrow
g> \frac{p_a+2}{2}.
\end{equation}

In particular, a consequence of Theorem \ref{thm:ijm} is
the following:

\begin{thm} \label{thm:eff}
Let $S$ be a smooth, projective $K3$ surface with $\Pic (S) \iso \ZZ[H]$ and $p_a:=p_a(H) \geq 2$. 
Let $Y$ and $\PP^1_{\Delta}$ be 
the generators of $N_1(\ss)_{\RR}$ with notation as in \S\;\ref{S:prel}.

If $X \in N_1(\ss)_{\ZZ}$ with 
$X \sim_{alg} Y -k\PP^1_{\Delta}$, then  $k \leq \frac{p_a+4}{4}$.
\end{thm}

\begin{proof}
We can assume that $X$ is an irreducible curve.
Then, precisely as in the case of a rational curve, $X$ corresponds either to the data of an 
irreducible curve $C \in |H|$ on $S$, with a partial normalization
$\widetilde{C}$ admitting a $2:1$ morphism onto the normalization $\widetilde{X}$ of $X$, or to the data
of an irreducible curve $C \in |H|$ on $S$ together with a point $x_0:=x_X \in S$. (The case corresponding to case (II)
in \S\;\ref{S:corr} does not occur, since the coefficient of $Y$ is one, precisely as in the case of a rational $X$
explained above.)

In the latter case $\mu(X) = \{{x_0}+C\} \subset \Sym^2(S)$, where $\mu: \ss \to \Sym^2(S)$ is the Hilbert-Chow morphism as usual, 
and one easily computes $k=(1/2)\mult_{x_0}(C)$ as in the rational case above. 
Since clearly $\mult_{x_0}(C) \leq 2$ if $p_a=2$ and $\mult_{x_0}(C) \leq 3$ if $p_a=3$, we have
 $k \leq \frac{p_a+4}{4}$ in these two cases. If $p_a \geq 4$, then from $\dim |H|-3-(p_a-4) =1$
and the fact that being singular at a given point imposes at most three independent conditions on 
$|H|$, we can find an irreducible curve $C' \in |H|$, different from $C$, singular at ${x_0}$, and passing through at least $p_a-4$ points of $C$. Therefore
\[ 2p_a-2=H^2=C'.C \geq \mult_{x_0}(C') \cdot \mult_{x_0}(C) +p_a-4 \geq 2\mult_{x_0}(C)+p_a-4,\]
whence $\mult_{x_0}(C) \leq (p_a+2)/2$, so that $k \leq (p_a+2)/4$.
 
In the first case, then, precisely as in the rational case above, 
\begin{equation} \label{eq:kappa}
k = \frac{p_a(\widetilde{C})+1}{2} -p_g(X)
\end{equation}
from Riemann-Hurwitz. 
By Brill-Noether theory on $\widetilde{X}$, it follows that $\widetilde{C}$ carries a $\mathfrak{g}^1_d$, with
\[ d \leq 2 \lfloor \frac{p_g(X)+3}{2} \rfloor. \]
By Theorem \ref{thm:ijm} we have $\rho_{sing}(p_a(C),1, d, p_a(\widetilde{C})) \geq 0$, whence $p_a(\widetilde{C}) \leq d-1 + p_a(C)/2$. 
The desired result now follows. 
\end{proof}

 By the proof of Theorem \ref{thm:eff} we see that if $C \in |mH|$ is an irreducible curve and $x_0 \in C$, then the class of the corresponding curve $\mu_*^{-1}\{{x_0}+C\} \subset \ss$ is given by $mY-(1/2)\mult_{x_0}(C) \PP^1_{\Delta}$. Hence
\[  \slope(\NE (\ss)) \leq \inf_{m \in \NN}\Big( \inf_{C \in |mH|} \Big( \inf_{x \in C} \frac{2m}{\mult_x(C)} \Big) \Big) =
   \inf_{m \in \NN} \frac{2}{H^2} \Big( \inf_{C \in |mH|} \Big( \inf_{x \in C} \frac{C.H}{\mult_x(C)} \Big) \Big). \]
It follows that
\begin{equation} \label{eq:seshadri}
\slope(\NE (\ss)) \leq \frac{\varepsilon(H)}{p_a-1}, 
\end{equation}
where 
\[ \varepsilon(H):= \inf_{x \in S} \Big( \inf_{C \ni x} \frac{C.H}{\mult_x(C)} \Big) \]
(and the infimum is taken over all irreducible curves $C \subset S$ passing through $x$)
 is the {\it (global) Seshadri constant} of $H$ (cf.
\cite[\S\;6]{dem}, \cite{ekl} or \cite{ba2}). These constants are very difficult to compute. The only case where they have been 
computed on general $K3$ surfaces is the case of quartic surfaces, where one has $\varepsilon(H)=2$
by \cite{ba1}, yielding the bound
$\slope(\NE (\ss)) \leq 1$. As a comparison, the bound one gets from \eqref{eq:slope} using the singular curves of genus two in $|H|$
is $\slope(\NE (\ss)) \leq 2/3$.
However, it is well-known that $\varepsilon(H) \leq \sqrt{H^2}$ on any surface, see e.g. \cite[Rem. 1]{ste}. Hence, by
\eqref{eq:seshadri} we obtain

\begin{thm} \label{thm:seshadribound}
Let $(S,H)$ be a primitively polarized $K3$ surface of genus $p_a:=p_a(H) \geq 2$ such that $\Pic(S) \iso \ZZ[H]$.
Then (cf. \eqref{inf}) 
 \begin{equation}\label{eq:seshadribound}
 \slope(\NE (\ss)) \leq \frac{\varepsilon(H)}{p_a-1} \leq \sqrt{\frac{2}{p_a-1}}.
\end{equation}
\end{thm}

In particular, \eqref{eq:seshadribound} shows that there is no lower bound on the slope of the Mori cone of $\ss$ of $K3$ surfaces, as the degree of the polarization tends to infinity, that is,
\begin{equation} \label{eq:inf0}
 \inf \Big\{ \slope(\NE (\ss)) \; | \; \mbox{$S$ is a projective $K3$ surface} \Big\} =0, 
\end{equation}
The same fact about $\slope_{rat}(\NE (\ss))$ will follow from \eqref{eq:slope-n} and
\eqref{eq:slope-d} below.

Note that one always has $\varepsilon(H) > \lfloor \sqrt{H^2} \rfloor-1$ under the hypotheses of Theorem \ref{thm:seshadribound}. Indeed, if $\varepsilon(H) < \sqrt{H^2}$, then there is an $x \in S$ and an irreducible curve $C$ such that $\varepsilon(H) = \frac{C.H}{\mult_x(C)}$, see e.g. \cite[Cor. 2]{og}. Since one easily computes $\dim |H \* \I_x^{(\lfloor \sqrt{H^2} \rfloor-1)}| \geq 2$, we can find a $D \in |L|$ such that 
$D \not \sup C$, $\mult_x (D) \geq \lfloor \sqrt{H^2} \rfloor-1$ and $D$ passes through at least one additional point of $C$. Thus
\[ \varepsilon(H) = \frac{C.H}{\mult_x(C)} = \frac{C.D}{\mult_x(C)} \geq 
\frac{\mult_x(C) \cdot \mult_x(D)+1}{\mult_x(C)} > \mult_x(D) \geq \lfloor \sqrt{H^2} \rfloor-1, \]
as desired. It follows that
 \begin{equation}\label{eq:sesbest}
 \frac{\varepsilon(H)}{p_a-1} > \frac{\lfloor \sqrt{2p_a-2} \rfloor -1}{p_a-1}, \; \mbox{for $(S,H)$ as in Theorem \ref{thm:seshadribound},}
\end{equation}
showing that there is a natural limit to how good a bound one can get on 
$\slope(\NE (\ss))$ by using Seshadri constants. 

The bound in \eqref{eq:seshadribound} is not (necessarily) obtained by rational curves in $\ss$. However, the presence of $p_g(X)$ in \eqref{eq:kappa}
above tends to indicate that the better bounds will be obtained by rational curves in $\ss$.
(Of course, if the Mori cone is closed, then the bound will indeed be obtained by rational curves, as explained at the end of \S\;\ref{S:prel}.)
In fact, the bound \eqref{eq:seshadribound} above will be improved, for infinitely many values of $H^2$, in 
Propositions \ref{prop:IP2} and \ref{prop:IP1bundle} below by rational curves. 

We now return to the study of irreducible {\it rational} curves in $\ss$ and to $\slope_{rat}(\NE (\ss))$.

Given Theorem \ref{thm:ijm} and \eqref{rho<0}, a natural question to ask is the following:

\medskip
\noindent
{\bf Hyperelliptic existence problem (HEP).} For $3 \leq p_g \leq \frac{p_a+2}{2}$,
does there exist a singular curve
in $|H|$ with hyperelliptic  normalization of geometric genus $p_g$?
\medskip

 By \eqref{eq:slope}  we have that
\begin{eqnarray} \label{eq:boundinf}
\textrm{a positive solution to (HEP) for ``maximal''} \; p_g= \lfloor \frac{p_a+2}{2} \rfloor \; \Longrightarrow \\
\nonumber \slope_{rat}(\NE (\ss)) \leq 
\begin{cases} 
\frac {4}{p_a+4}  & \; \mbox{if $p_a$ is even;} \\ 
\frac {4}{p_a+3}  & \; \mbox{if $p_a$ is odd} 
\end{cases} 
\end{eqnarray}
and, by \eqref{eq:square}, the $q$-square of the associated rational curves would be much less
than what predicted by Hassett and Tschinkel \cite[Conj. 3.1]{HT}.  Moreover, the bounds in 
\eqref{eq:boundinf} would be much stronger than the bound given by the right hand inequality  in
\eqref{eq:seshadribound}, and even stronger than the best bounds one could obtain from Seshadri constants
(compare the left hand side inequality in \eqref{eq:seshadribound} with \eqref{eq:sesbest}).

It is natural to try to solve
(HEP) using nodal curves, as one has better control of their deformations and
their parameter 
spaces (the Severi varieties considered in \S\;\ref{S:esistenza3}). 
 After the positive answer to the hyperelliptic existence problem
for the specific values $p_g=3$ and $p_a=4,5$ in
\cite[Examples 2.8 and 2.10]{fkp}, Theorem \ref{thm:esistenza} gives the first examples, 
at least as far as we know, of positive answers to 
the hyperelliptic existence problem for primitively polarized $K3$ surfaces of {\it any} degree.  

In Remark \ref{rem:g=g0} we showed that  $p_g(C)=g_0(C) =3$ for these constructed curves $C \in |H|$ (cf. \eqref{eq:defg0}), so that the classes of the associated rational curves $R_C \subset \ss$ are, 
using \eqref{eq:w_X}, 
\begin{equation} \label{eq:genere3}
w_{R_C} = H-2\mathfrak{e}, 
\end{equation}
with
\[ q(w_{R_C})=q(R_C)=2p-10 \geq -2. \]
Moreover, using \eqref{eq:slope}, Theorem \ref{thm:esistenza} yields (cf. \eqref{inf-rat}):

\begin{cor} \label{cor:g=3}
Let $(S,H)$ be a general, primitively polarized $K3$ surface of genus $p_a(H) \geq 4$.
Then 
 \begin{equation} \label{eq:g=3boundinf}
 \slope_{rat}(\NE (\ss)) \leq \frac {1}{2}.
\end{equation}
\end{cor}

Note that the existence of nodal curves of geometric genus $2$ in $|H|$, which was already known and followed from the  nonemptiness of the Severi varieties 
on general $K3$ surfaces, as explained in the beginning of 
\S\;\ref{S:esistenza3},  leads to the less good bound of $\frac{2}{3}$.  Therefore, again as far as we know, \eqref{eq:g=3boundinf} is the first ``nontrivial'' bound on 
the slope of rational curves holding for all degrees of the polarization.
 As already mentioned, for infinitely many degrees of the polarization we will in fact improve this bound in Propositions \ref{prop:IP2} and \ref{prop:IP1bundle} below. 

\begin{remark}\label{rmk:mH} 
{\rm One may also look for irreducible singular curves with hyperelliptic
normalizations
 in $|mH|,\ m\geq 2$.  In \cite[Corollary 4]{fkp}, we also proved that,
 apart from some special numerical cases (where we were not able to
conclude),
 the negativity of $\rho_{sing}(p_a(mH),1,2,g)$ 
 implies the
 non-existence of {\it irreducible nodal}  curves in $|mH|$ with hyperelliptic normalizations.
 A positive solution to the hyperelliptic existence problem for singular
curves in $|mH|$  would then provide an even better bound on the slope of
the Mori cone.
 Namely, one would for instance get 
$ \slope(\NE (\ss))  \leq {4}/[{m(p_a(H)+4)}]$ for even $p_a$.
Whereas we tend to believe that the nonnegativity of $\rho_{sing}$
should imply existence of curves with hyperelliptic normalizations for the specific values of $p_a$ and $g$ in a primitive linear system $|H|$ on a general $K3$, we are not  sure what to expect 
 for curves in $|mH|$ when $m >1$. For instance, the degeneration methods to prove existence as in the proof of Theorem 
\ref{thm:esistenza} will certainly get more difficult, because the irreducibility of the obtained curves after deformation is not automatically ensured.}
\end{remark}

\begin{remark}\label{rmk:deform} 
{\rm We do not know whether there will always be components in $|H|^{hyper}$ (whenever nonempty) of singular curves with hyperelliptic normalizations
such that the singularities of the general member are as nice as possible, that is, all nodes and all non-neutral \cite[\S 3]{fkp}.}
\end{remark}

%
\section{$\PP^2$s and threefolds birational to $\PP^1$-bundles  in the Hilbert square of a general $K3$ surface}\label{S:esempi}
%

We now give an infinite series of examples of {\it general}, primitively polarized $K3$ surfaces
$(S,H)$, of infinitely many degrees such that $\ss$ 
contains either a $\PP^2$ or a threefold birational to a $\PP^1$-bundle, thus showing both possibilities occurring in Proposition \ref{cl:dimRV}. 

Both series of examples are similar to Voisin's constructions in
\cite[\S\;3]{V2}. The idea is to start with a smooth quartic surface $S_0$ such that $S_0^{[2]}$ contains an ``obvious''
$\PP^2$ or threefold birational to a $\PP^1$-bundle over $S_0$, use the involution on the quartic to produce another such
$\PP^2$ or uniruled threefold, and then deform $S_0$ keeping the latter one and loosing the first one in the Hilbert square.

We remark that the question of existence of $\PP^2$s in $\ss$ when $S$ is $K3$ is a very interesting problem because of the following fact: a $\PP^2$ in $\ss$ gives rise to a birational map from $\ss$ onto another hyperk\"ahler
fourfold, and conversely any birational transformation
$X --\khpil X'$ between projective, symplectic fourfolds can be factorized into a finite sequence of Mukai flops (cf. \cite[Thm. 0.7]{Muk}), by \cite[Thm. 2]{W2}, see also \cite{BHL,HY,WW}. 
Therefore, in the case of a $K3$ surface, if $\ss$ contains no $\PP^2$s, then $\ss$ admits no other birational model than itself. 

Also uniruled divisors have an influence on the  birational geometry of a hyperk\"ahler manifold $X$. Indeed, Huybrechts proved in \cite[Prop. 4.2]{H2} that a class $\alpha$
in the closure of the positive cone 
$\overline{\mathcal C}_X$ lies in the closure of the 
birational K\"ahler cone $\overline {\mathcal {BK}}_X$
if and only if $q(\alpha,D)\geq0$, for all uniruled divisors $D\subset X$. (Recall that the {\it positive cone}
$\mathcal{C}_X$ is the connected component of 
$\{ \alpha\in H^{1,1}(X,\RR) : q(\alpha)\geq 0\}$
containing the cone $\mathcal{K}_{X}$ of all K\"ahler classes of $X$, and the {\it birational K\"ahler cone} $\overline {\mathcal {BK}}_X$ equals by definition
$\cup_{f:X- - \to X'} f^*\mathcal{K}_{X'}$,  
where $f$ is a bimeromorphic map onto another hyperk\"ahler manifold $X'$).

%
\subsection{$\PP^2$s in $\ss$} \label{S:hassett}
%
The first nontrivial case, the case of degree $10$, is particularly easy, so we begin with that one. 

\begin{exa} \label{exa:d=10}
{\rm (Hassett) Let $S\subset \mathbb P^6$ be a general 
$K3$ surface of degree $10$. By \cite{mukk3} the surface $S$ is a complete intersection $S= G \cap T \cap Q$,
where $G:=\Grass(2,5)$ is the Grassmannian of lines in $\mathbb P^4$ embedded in $\mathbb P^9$ by
its Pl\"ucker embedding, $T$ is a general $6$-dimensional linear subspace of $\mathbb P^9$,
and $Q$ is a hyperquadric in $\mathbb P^9$.
Set $Y:=G\cap T$. Then $Y$ is a Fano $3$-fold of index $2$. Let $F(Y)$ be its variety of lines.
It is classically known (see e.g. \cite{Fa} for a modern proof) that $F(Y) \cong \mathbb P^2$.
 Then we may embed this plane in $S^{[2]}$ by mapping the point corresponding to  a line $[\ell]$ 
to $\ell \cap Q$. By generality, $S$ does not contain any line, so
that this map is a morphism.}
\end{exa}

\vspace{0,3cm}
The construction behind the following result, generalizing the previous example, was shown to us by B.~Hassett.

\begin{prop}\label{prop:IP2}
Let $(S,H)$ be a general primitively polarized $K3$ surface of degree $H^2=  2(n^2-9n+19)$, for $n \geq 6$.
Then $S^{[2]}$ contains a $\PP^2$. 

The class $w_{\ell} \in H^2(\ss, \QQ)$ corresponding to a line $\ell\subset \PP^2$ is
\begin{equation} \label{eq:rettaq}
w_{\ell} = H - \frac{2n-9}{2} \mathfrak{e},
\end{equation}
In particular
\begin{equation} \label{eq:slope-n}
 \slope_{rat}(\NE (\ss)) \leq \frac{2}{2n-9}. 
\end{equation}

Moreover the curves $C\subset S$ with hyperelliptic normalizations 
associated to the lines $\ell\subset \PP^2\subset \ss$ lie in $|H|$, have 
geometric genus $p_g=2n-10$, and $\rho_{sing}(p_a(C),1,2,p_g)=n(n-13)+42\geq 0$.

\end{prop}

\begin{proof}  Consider the lattice $\ZZ F \+ \ZZ G$ with
intersection matrix
\[  \left[
  \begin{array}{cc}
   F^2        &   F.G                   \\
   G.F        &   G^2
    \end{array} \right]  =
    \left[
  \begin{array}{cc}
  2 & n       \\
  n &  4
    \end{array} \right], \; n \geq 6.     \]
Since it has signature $(1,1)$, then, by a result of Nikulin \cite{nik}
(see also \cite[Cor. 2.9(i)]{Morr}), there is an algebraic $K3$
surface $S_0$ with the given Picard lattice. Performing
Picard-Lefschetz reflections on the lattice, we can assume that $G$ is nef,
by \cite[VIII, Prop. 3.9]{BPV}. By Riemann-Roch and Serre duality,
we have $G >0$ and $F>0$.
Straightforward computations on the Picard lattice rules out the existence of divisors $\Gamma$ satisfying $\Gamma^2=-2$ and $\Gamma.F \leq 0$ or $\Gamma.G \leq 1$;  or $\Gamma^2=0$ and $\Gamma.F=1$ or $\Gamma.G=1,2$. By \cite{SD} it follows that both $|F|$ and $|G|$ are base point free, 
$\varphi_{|F|}: S_0 \to \mathbb P^2$
is a double cover and $\varphi_{|G|}: S_0 \to \mathbb P^3$ is an embedding onto a smooth quartic not containing lines. As explained in \S\;\ref{S:bound}, $S_0^{[2]}$ contains
a $\mathbb P^2$ arising from the double cover.

If $\ell_0$ is a line on the $\mathbb P^2$, the corresponding class  in
$H^2(S_0^{[2]}, \mathbb Q)$ is $w_{\ell_0}= 2F-3\mathfrak{e}$, which coincides with the corresponding integral class 
$\rho_{\ell_0}$ (cf. \cite[Example 5.1]{HT}).

As $S_0$ is a quartic surface not containg lines, $S_0^{[2]}$ admits an
involution
\[
 \iota :S_0^{[2]}\to S_0^{[2]}; \; \; \xi \mapsto (\ell_\xi \cap S_0) \setminus \xi,
\]
by \cite[Prop. 11]{B2},
where $\ell_\xi$ is the line determined by $\xi$, and the sign $\setminus$ means that we take the
residual subscheme. The corresponding involution on cohomology is given by (cf. e.g. \cite[(4.1.6)-(4.1.7)]{ogr})
\[
 v\mapsto q(G-\mathfrak{e},v)\cdot (G-\mathfrak{e}) - v.
\]
The involution sends the $\mathbb P^2$ into another $\mathbb P^2$, and the corresponding class
associated to a line on it is
\begin{equation} \label{eq:classeretta}
 q(G-\mathfrak{e},2F-3\mathfrak{e})\cdot(G-\mathfrak{e}) - (2F-3\mathfrak{e}) = 2((n-3)G-F) - (2n-9)\mathfrak{e}.
\end{equation}
In order to obtain a general $K3$ with the desired property we now deform $S_0^{[2]}$.
Precisely, we consider a general deformation of $S_0^{[2]}$ such that
(i) $\mathfrak{e}$ remains algebraic
and (ii) $\iota(\mathbb P^2)$ is preserved. Deformations satisfying (i) form
a countable union of hyperplanes in the deformation space of
$S_0^{[2]}$, which is smooth and of dimension $21$, and may be characterized as those
of the form $S^{[2]}$, where $S$ is a $K3$ surface (see \cite[Thm. 6 and Rem. 2]{B}). Deformations preserving $\iota(\mathbb P^2)$ can be characterized as those preserving the image in $H^2(S^{[2]}, \ZZ)$ of the class of the line in
$\iota(\mathbb P^2)$ as an algebraic class (see \cite[Thm. 4.1 and Cor. 4.2]{HT} or \cite{V2}), that is, using \eqref{eq:classeretta}, those deformations keeping
$H:=(n-3)G-F \in \Pic (S_0^{[2]})$, or, equivalently, $H \in \Pic (S)$, by
\eqref{eq:decpic}. As $H^2 = [(n-3)G-F]^2= 2(n^2-9n+19) \geq 2$ for $n \geq 6$  and $H$ is primitive, those deformations form a divisor in the $20$-dimensional space
of deformations keeping $\mathfrak{e}$ algebraic, by \cite[Thm. 14]{Kod}.

We therefore obtain a $19$-dimensional space of deformations of $S_0^{[2]}$,
whose general member is $S^{[2]}$, where $(S,H)$ is a general primitively polarized
(algebraic) $K3$ surface of degree $H^2=2(n^2-9n+19)$, $n \geq 6$, and $S^{[2]}$ contains a plane.

The class $w_{\ell} \in H^2(\ss, \QQ)$ corresponding to the line $\ell$ is as in 
\eqref{eq:rettaq}, yielding \eqref{eq:slope-n}.

As $S$ is general, it does not contain
smooth rational curves, so that the $\PP^2$ is not of the form $C^{[2]}$, for a smooth rational curve $C$ on $S$. By Lemma \ref{lem:famcost}, the lines in the $\PP^2$ in $S^{[2]}$ give rise to a two-dimensional family $V$ of curves on $S$ with hyperelliptic normalizations,
so that $R_V = \mu(\PP^2)$, where $\mu: S^{[2]} \to \Sym^2(S)$ is the Hilbert-Chow morphism. 
By \eqref{eq:rettaq} we have $\ell.H=H^2$, so that, by the very definition of the divisor $H$ in $H^2(\ss,\ZZ)$, the lines in the $\PP^2$ correspond to curves $C \in |H|$. Comparing \eqref{eq:w_X}
 and \eqref{eq:rettaq}, we see that $g_0(C)=2n-10$, cf. \eqref{eq:defg0}. Now we note that the 
general line in the $\PP^2$ is not tangent to $\Delta=2\mathfrak{e}$. (Indeed, this follows 
by deformation since in $\ss_0$ we have that $\iota(\PP^2) \cap \Delta$ is a smooth plane sextic,
since we have a composite map $S_0 \khpil \PP^2 \khpil \iota(\PP^2)$ that is finite of degree two, whence ramified along a smooth sextic, as $S_0$ is a smooth $K3$.)
Therefore we have $p_g(C)=2n-10$. We compute
$\rho_{sing}=n(n-13)+42 \geq 0$ (recall that $n \geq 6$). 
\end{proof}

The examples contained in the above proposition is interesting in several regards.

Notice first that $q(\ell)=-5/2$, cf. \eqref{eq:qalfa}, in accordance with 
the prediction in \cite[Conj. 3.6]{HT}.
 
The proposition shows in particular that the correspondence in Remark \ref{rem:apppg2} is not one-to-one
 and also shows that the case $\dim (V)=\dim(R_V)=2$ of Proposition \ref{cl:dimRV} actually occurs.
 
The result also gives nontrivial examples of curves in $|H|$ with hyperelliptic normalizations and positively answers 
the hyperelliptic existence problem for
$p_a= n^2-9n+20$ and $p_g=2n-10$, $n \geq 6$.

Moreover \eqref{eq:slope-n} shows that there is no lower bound on 
$\slope_{rat}(\NE (\ss))$ as the degree of the polarization tends to infinity. 
The same follows from 
\eqref{eq:slope-d} in Proposition \ref{prop:IP1bundle} below. Both the bounds  \eqref{eq:slope-n} and \eqref{eq:slope-d} below
in fact yield better bounds on $\slope(\NE (\ss))$ than \eqref{eq:seshadribound}.

Finally, the
conics on the $\PP^2$ give a five-dimensional family $V(2)$ of irreducible curves with hyperelliptic normalizations on $S$. Of course this family has obvious non-integral members, corresponding to non-integral conics. More generally, for any $m \geq 3$, the $(3m-1)$-dimensional family of nodal rational curves in $|\cO_{\PP^2}(m)|$ (cf. \cite[Thm. 1.1]{CS}) yields corresponding
families $V(m)$ of curves in $|mH|$ with hyperelliptic normalizations
with $\dim V(m) =3m-1 \geq 5$ and $\dim (R_V)=2$, showing in particular that
the case $\dim (V) >\dim(R_V)=2$ of Proposition \ref{cl:dimRV} actually occurs.

In the case of the conics, we compute 
$p_g=4n-19$ as above and as $p_a(2H)=4n^2-36n+77$, we
get $\rho_{sing}=4n(n-11)+117 \geq -3$ in these cases. This does not contradict 
\cite[Thm. 1]{fkp}. 

%
\subsection{Threefolds birational to $\PP^1$-bundles in $\ss$} \label{S:andreas}
%

We start with an explicit example in the special case of a quartic surface.

\begin{exa} \label{esempigianluca4}{\rm In the case of a general quartic $S$ in $\PP^3$ we can find a $\PP^1$-bundle over $S$ in $\ss$, arising from the two-dimensional family of hyperplane sections of geometric genus two. In fact, taking the tangent plane
through the general point of $S$ we get a nodal curve of geometric
genus $2$. We obtain in this way a family $V$ of nodal
curves with hyperelliptic normalizations in the hyperplane linear system. This family
is parametrized by an open subset of $S$, and the locus in $\ss$ covered by the associated rational curves is  birational to a $\P^1$-bundle over this open
subset. To see this, set $C_p:= (S \cap T_p S )$, and let
$\widetilde{C}_p$ be the normalization of $C_p$.
Note that the $\mathfrak{g}^1_2$ on $\widetilde{C}_p$, viewed on $C_p$,
is given by the pencil of lines in $T_p S$ through the
node $p$. If, for two distinct points $p,q \in S$, the
$\mathfrak{g}^1_2$s on $\widetilde{C}_p$ and $\widetilde{C}_q$ had
two common points, say $x$ and $y$ (so that the map $\Phi_V$
in \eqref{eq:univ2} sends $(p,x+y)$ and $(q,x+y)$ to the same point $x+y$ in
$\Sym^2(S)$),
then  the line $T_p S \cap T_q S$, which is bitangent to $S$, 
would also pass through $x$ and $y$. This is absurd, as $\deg(S)=4$.

By \eqref{eq:w_X}, the class $w \in H^2(\ss, \QQ)$ corresponding
to the curves of geometric genus $2$ is $w= H-\frac{3}{2}\mathfrak{e}$, whence 
$q(w)=-1/2$, 
as predicted by \cite[Conj. 3.6]{HT}. 
Moreover, performing the usual involution on the quartic, we send
the constructed uniruled threefold to another one, with corresponding fibre class given by $\mathfrak{e}$, so that it simply is the $\PP^1$-bundle $\Delta$ over $S$. 
This shows that also our original threefold was smooth, so in fact a $\PP^1$-bundle over $S$.
}
\end{exa}

\vspace{0,3cm}

We now give  an infinite series of examples of general $K3$s whose Hilbert squares contain threefolds birational to $\PP^1$-bundles.

\begin{prop}\label{prop:IP1bundle}
Let $(S,H)$ be a general primitively polarized $K3$ surface of degree $H^2=  2(d^2-1)$, for $d \geq 2$.
Then $S^{[2]}$ contains a threefold birational to a $\PP^1$-bundle over a $K3$ surface. 

The class $w_{f} \in H^2(\ss, \QQ)$ corresponding to a fibre  is
\begin{equation} \label{eq:fibraz}
w_{f} = H - {d} \mathfrak{e}\in  H^2(\ss, \ZZ).
\end{equation}
In particular
\begin{equation} \label{eq:slope-d}
 \slope_{rat}(\NE (\ss)) \leq \frac{1}{d}. 
\end{equation}

Moreover the curves $C\subset S$ with hyperelliptic normalizations 
associated to the fibres of the threefold lie in $|H|$, have 
geometric genus $p_g=2d-1$, and $\rho_{sing}(p_a(C),1,2,p_g)=d(d-4)+4\geq 0$.

\end{prop}

\begin{proof}
This time we start with the lattice $\ZZ F \+ \ZZ G$ with
intersection matrix
\[  \left[
  \begin{array}{cc}
   F^2        &   F.G                   \\
   G.F        &   G^2
    \end{array} \right]  =
    \left[
  \begin{array}{cc}
  -2 & d      \\
  d &  4
    \end{array} \right], \; d \geq 2.     \]
As in Proposition \ref{prop:IP2} one easily shows that there is an algebraic $K3$
surface $S_0$ with $\Pic (S_0) = \ZZ F \+ \ZZ G$ and that
$\varphi_{|G|}: S_0 \to \mathbb P^3$ is an embedding onto a smooth quartic not containing lines and $F$ is a smooth, irreducible rational curve.
(Note that $F^{[2]} = \PP^2$ and performing the same procedure  on this plane as in the proof of Proposition \ref{prop:IP2}, one gets precisely the same series of examples as above.)

We now consider the divisor $F \subset S_0^{[2]}$, defined as the length-two schemes with some support along $F$. One easily sees that this is a threefold birational to a $\PP^1$-bundle over $S_0$ and that the class in $H^2(S_0^{[2]},\ZZ)$ corresponding to the fibres $f$ is $\rho_f=F$, cf. \cite[Example 4.6]{HT}.

The involution on the quartic sends this threefold to another threefold birational to a $\PP^1$-bundle over $S_0$ and the corresponding class of the fibres is
\begin{equation} \label{eq:classefibra}
 q(G-\mathfrak{e},F)\cdot(G-\mathfrak{e}) - F = dG-F-d\mathfrak{e}.
\end{equation}
Note that this threefold satisfies the conditions in \cite[Thm. 4.1]{HT} by \cite[Example 4.6]{HT}, so that, as in the previous example, we can deform $S_0^{[2]}$, keeping $\mathfrak{e}$ algebraic and $H:=dG-F$. We thus obtain a $19$-dimensional space of deformations of $S_0^{[2]}$,
whose general member is $S^{[2]}$, where $(S,H)$ is a general, primitively polarized
(algebraic) $K3$ surface of degree $H^2=2(d^2-1) \geq 6$ and
$S^{[2]}$ contains a threefold birational to a $\PP^1$-bundle, again over a $K3$ surface (see also \cite[Thm. 4.3]{HT}).

The unique class $w_f \in H^2(\ss, \QQ)$ corresponding to a fibre $f$ is as in 
\eqref{eq:fibraz} and yields \eqref{eq:slope-d}.

By \eqref{eq:fibraz} we have $f.H=H^2$, so that, by the very definition of the divisor $H$ in $H^2(\ss,\ZZ)$, 
the fibres $f$
of $Y$ correspond to curves $C \in |H|$. Comparing \eqref{eq:w_X} and 
\eqref{eq:fibraz}, we see that $g_0(C)=2d-1 \geq 3$, cf. \eqref{eq:defg0}. As in the proof of Proposition \ref{prop:IP2}, one can see that the general fibre of $Y$  is not tangent to $\Delta=2\mathfrak{e}$, so that in fact we have $p_g(C)=2d-1$. In particular, $Y$ is not one of the obvious uniruled threefolds arising from the rational curves on $S$, or the one-dimensional families of elliptic curves on $S$. A computation shows that
$\rho_{sing}=d(d-4)+4 \geq 0$.
\end{proof}

Again, a few comments are in order.

The square of the class of the fibres of the uniruled threefolds constructed above 
is $q(f)=-2$,  as predicted in \cite[Conj. 3.6]{HT}.

The obtained family $V$ of curves on $S$ with hyperelliptic normalizations
has  $\dim (V)=2$ and $\dim (R_V)=3$, showing that also this case of Proposition \ref{cl:dimRV} actually occurs.
This family gives nontrivial examples of curves in $|H|$ with hyperelliptic normalizations and positively answers the hyperelliptic existence problem for
$p_a=2(d^2-1)$ and  $p_g=2d-1$ for every $d \geq 2$.
Note that the case $d=2$ is the case described in \cite[Example 2.8]{fkp}.

%
%

\vskip 30pt

\noindent
{\small Flaminio Flamini, Dipartimento di Matematica, Universit\`a degli Studi di Roma "Tor Vergata", Viale
della Ricerca Scientifica, 00133 Roma, Italy. e-mail {\tt flamini@mat.uniroma2.it.}}

\vskip 10pt

\noindent
{\small Andreas Leopold Knutsen, Dipartimento di Matematica, 
Universit\`a Roma Tre, Largo San 
Leonardo Mu\-rial\-do 1, 00146, Roma, Italy. 
e-mail {\tt knutsen@mat.uniroma3.it.}}

\vskip 10pt

\noindent
{\small Gianluca Pacienza, Institut de Recherche Math\'ematique Avanc\'ee,
Universit\'e L. Pasteur et CNRS, rue R. Descartes - 67084 Strasbourg Cedex, France. e-mail {\tt pacienza@math.u-strasbg.fr.}}

\newpage

\pagestyle{fancy}

\appendix

\rhead{{\scriptsize \thepage}}
\chead{{\scriptsize EDOARDO SERNESI: PARTIAL DESINGULARIZATIONS OF FAMILIES OF NODAL CURVES}}
\lhead{}

\lfoot{}
\cfoot{}
\rfoot{}

%
\section{} \label{A:edoardo}
%

\begin{center}
{\bf PARTIAL DESINGULARIZATIONS OF FAMILIES OF NODAL CURVES}   \\
\vspace{5mm}
{\rm {\small EDOARDO SERNESI}}\footnote{Work done during a visit to the Institut Mittag-Leffler (Djursholm, Sweden), whose support is gratefully acknowledged. I am grateful to F.~Flamini, A.~L.~Knutsen and G.~Pacienza for accepting 
this note as an Appendix to their paper, and to F.~Flamini for some useful remarks.} 
\end{center}

\vskip 10pt

In this Appendix we show how to construct simultaneous partial desingularizations of families of nodal curves, generalizing a well known procedure of simultaneous total desingularization, as described in \cite{bT80}.  

 We work over an algebraically closed field {\bf k} of characteristic $0$. For every morphism $X \to Y$, and for every $y \in Y$, we denote by  $X(y)$ the scheme-theoretic fibre of $y$.  

\begin{theorem}\label{thm:sern}
 Let 
\[
 \xymatrix{
f: \C \ar[r] & V
}
\]
 be a flat projective  family of curves, with $\C$ and $V$ algebraic schemes, such that all fibres have at most ordinary double points (nodes) as singularities. Let $\delta \ge 1$ be an integer. Then there is a commutative diagram:
\[
\xymatrix{
    D_{\delta} \hspace{0,2cm} \ar@{^{(}->}[r]       \ar[dr]_{q} &  \C' \ar[d]^{f'} \ar[r] & \C \ar[d]^{f} \\
    & E_{(\delta)} \ar[r]^{\alpha} & V
}
\]
 with the following properties:
 
\begin{itemize}
\item[(i)] $\alpha$ is finite and unramified,   the square is cartesian, and $q$ is an \'etale cover of degree $\delta$.  
\item[(ii)]   The left triangle defines a marking of all $\delta$-tuples of nodes of fibres of $f$. 
 In particular  $f'$ parametrizes all curves of the family $f$ having $\ge \delta$ nodes
 and, for each $\eta \in E_{(\delta)}$, 
 $D_\delta(\eta)\subset \C'(\eta)$  is a set of $\delta$  nodes of the curve $\C'(\eta)$.
 \item[(iii)]  The diagram is universal with respect to properties (i) and (ii). Precisely, if  
 \[
\xymatrix{
    \tilde{D}  \hspace{0,2cm} \ar@{^{(}->}[r]       \ar[dr]_{\tilde{q}} &  \tilde{E} \x_V \C \ar[d]^{\tilde{f}} \ar[r] & \C \ar[d]^{f} \\
    & \tilde{E} \ar[r] & V
}
\]
 is a diagram having the properties analogous to (i) and (ii), then there is a unique factorization \[
\xymatrix{
\tilde E \ar[r]^{\varphi} & E_{(\delta)}  \ar[r]^{\alpha} & V  
}
\]
 such that $\tilde q$ and $\tilde f$ are obtained by pulling back $q$ and $f'$ by  $\varphi$.  
 \end{itemize}

 If moreover $E_{(\delta)}$ is normal, then the above diagram can be enlarged as follows:

\[
\xymatrix{
& \overline{\C} \ar[d]^{\beta} & \\ 
    D_{\delta}  \hspace{0,2cm} \ar@{^{(}->}[r]       \ar[dr]_{q} &  \C' \ar[d]^{f'} \ar[r] & \C \ar[d]^{f} \\
    & E_{(\delta)} \ar[r]^{\alpha} & V
}
\]
 where:
 
 \begin{itemize}
 
 \item[(iv)] $\beta$ is a birational morphism such that, for each $\eta \in E_{(\delta)}$, the restriction:
 \[
 \xymatrix{
\beta(\eta): \overline{\C}(\eta) \ar[r] & \C'(\eta)
}
\]
 is the partial normalization at the nodes $D_\delta(\eta)$.

 \item[(v)]  The composition $\bar f := f'\circ\beta$ is flat.
 
 \end{itemize}
 
\end{theorem}

\lhead{{\scriptsize \thepage}}
\chead{{\scriptsize EDOARDO SERNESI: PARTIAL DESINGULARIZATIONS OF FAMILIES OF NODAL CURVES}}
\rhead{}

\lfoot{}
\cfoot{}
\rfoot{}

\begin{proof} 
 Consider the first relative cotangent sheaf 
 $\T^1_{\C/V}$. Since all fibres of $f$ are nodal,  $\T^1_{\C/V}$ commutes with base change 
(\cite[Lemma 4.7.5]{eS06} or \cite{jW74}), thus on every fibre $\C(v)$, $v\in V$,  
 it restricts to  $\T^1_{\C(v)}$, which is the structure sheaf of the scheme of nodes of $\C(v)$.  It follows that we have 
 \[\T^1_{\C/V} = \Oc_E
 \]
 for a    closed subscheme $E \subset \C$ supported on the nodes of the fibres of $f$. Consider the composition   
 \[
 \xymatrix{
f_E: E \subset \C \ar[r]^{\hspace{0,7cm} f} & V
}
\]
    By construction it follows  that $f_E$ is finite and unramified. 
 Now fix $\delta \ge 1$ and consider the fibre product:
 \[
 \underbrace{E \times_V \cdots\times_VE}_\delta
 \]
 Since $f_E$ is finite and unramified,     it follows from \cite[Exp.1, Prop. 3.1]{SGA1}, and by induction on 
 $\delta$ (see \cite[Lemma 4.7.11(i)]{eS06}), that  we have a disjoint union decomposition:
 \[
 E \times_V \cdots\times_VE =  \Delta \coprod E_\delta
 \]
 where $\Delta$ is the union of all the diagonals, and $E_\delta$ consists of all the ordered $\delta$-tuples of distinct points of $E$ mapping to the same point of $V$; moreover the natural projection morphism 
\[
 \xymatrix{
E_\delta \ar[r] &  V
}
\]
 is finite and unramified. 
 
 There is a natural action of the symmetric group $\Sigma_\delta$ on $E_\delta$ that commutes  with the    projection to $V$.
 We denote the quotient 
 $E_\delta/\Sigma_\delta$ by $E_{(\delta)}$. Since the composition
\[
 \xymatrix{
E_\delta \ar[r] &  E_{(\delta)}  \ar[r] & V
}
\]
 is finite and unramified and the first morphism is an \'etale cover,  the morphism $ \alpha:E_{(\delta)} \to V$ is 
  finite and unramified.  
 Note that if, for  a closed point $v \in V$,  $\C(v)$ has $\delta+t$ nodes as the only singularities, 
 with $t >  0$, then $\alpha^{-1}(v)$ has degree ${\delta+t\choose t}$. 
 Now let 
 \[
 D_\delta = \{(\eta, e):   e\in \Supp(\eta)\} \subset E_{(\delta)}\times_V E
 \]
 Then the first projection defines the tautological family:
\begin{equation}\label{E:dia1}
 \xymatrix{ D_\delta \ar[d]^q & \hspace{-0,8cm} \subset \hspace{0,2cm}  E_{(\delta)}\times_V E  \hspace{0,2cm} 
\subset  \hspace{0,2cm}  E_{(\delta)}\times_V \C \\
E_{(\delta)} &
} 
 \end{equation}  
 which is an \'etale cover of degree $\delta$. The fibre $D_\delta(\eta)$ is the $\delta$-tuple parametrized by $\eta$, 
 for each 
 $\eta \in E_{(\delta)}$\footnote{If $\delta=1$ then  $E_{(1)}=E$ and $D_1  \subset E\times_VE$ is the diagonal.}. 
 We therefore have the following diagram:
 \[
\xymatrix{
    D_{\delta}  \hspace{0,2cm} \ar@{^{(}->}[r]       \ar[dr]_{q} &  \C' \ar[d]^{f'} \ar[r] & \C \ar[d]^{f} \\
    & E_{(\delta)} \ar[r]^{\alpha} & V
}
\]
 where we have denoted by $\C' = E_{(\delta)}\times_V\C$. 
The fibres of $f'$ are all the curves  of the family $f$ having  $\ge \delta$ nodes. 
For each $\eta\in E_{(\delta)}$ the divisor $D_\delta(\eta)\subset \C'(\eta)$ marks the set of $\delta$ nodes parametrized by $\eta$.  This proves (i) and (ii).

(iii) follows from the fact that  $\alpha: E_{(\delta)} \to V$ is the relative Hilbert scheme of degree $\delta$ of 
$f_E: E \to V$, and  (\ref{E:dia1}) is the universal family.

Assume that $E_{(\delta)}$ is normal. Then we can normalize $\C'$ locally around $D_\delta$ as in 
\cite[Theorem 1.3.2]{bT80},   to obtain a birational morphism $\beta$ having the required properties (iv) and (v).
 \end{proof}

  A typical example of the situation considered in the theorem is when $V$ parametrizes a complete linear system of curves on an algebraic surface. If the morphism $f_E$ is \emph{self-transverse of codimension 1} (see \cite[Definition 4.7.13]{eS06}) then  the Severi variety of irreducible $\delta$-nodal curves is nonsingular and of codimension $\delta$, and $E_{(\delta)}$ is nonsingular (see 
\cite[Lemma 4.7.14]{eS06}), so that the theorem applies and the simultaneous  partial desingularization exists.
  This happens for example for the linear systems of plane curves \cite[Proposition 4.7.17]{eS06}.  
 
 \rhead{{\scriptsize \thepage}}
\chead{{\scriptsize EDOARDO SERNESI: PARTIAL DESINGULARIZATIONS OF FAMILIES OF NODAL CURVES}}
\lhead{}

\lfoot{}
\cfoot{}
\rfoot{}

\vskip 30pt

\noindent
{\small Edoardo Sernesi, Dipartimento di Matematica, 
Universit\`a di Roma Tre, Largo San 
Leonardo Murialdo 1, 00146, Roma, Italy. 
e-mail {\tt sernesi@mat.uniroma3.it.}}

\end{document}